\documentclass[reqno,11pt]{amsart}
\usepackage{amsmath, latexsym, amsfonts, amssymb, amsthm, amscd}
\usepackage{graphics,epsf,psfrag}
\numberwithin{equation}{section}

\newtheorem{theorem}{Theorem}[section]

\newtheorem{rem}[theorem]{Remark}

\newcommand{\ind}{\mathbf{1}}

\newcommand{\R}{\mathbb{R}}
\newcommand{\Z}{\mathbb{Z}}
\newcommand{\N}{\mathbb{N}}
\renewcommand{\tilde}{\widetilde}

\newcommand{\med}[1]{\left\langle #1\right\rangle}

\newcommand{\cL}{{\ensuremath{\mathcal L}} }

\newcommand{\cD}{{\ensuremath{\mathcal D}} }

\newcommand{\bP}{{\ensuremath{\mathbf P}} }
\newcommand{\bE}{{\ensuremath{\mathbf E}} }

\DeclareMathSymbol{\leqslant}{\mathalpha}{AMSa}{"36} 
\DeclareMathSymbol{\geqslant}{\mathalpha}{AMSa}{"3E} 
\DeclareMathSymbol{\eset}{\mathalpha}{AMSb}{"3F}     
\newcommand{\dd}{\,\text{\rm d}}             
\newcommand{\limtwo}[2]{\lim_{\substack{#1 \\ #2}}}     

\newcommand{\bbE}{{\ensuremath{\mathbb E}} }

\newcommand{\bbP}{{\ensuremath{\mathbb P}} }

\newcommand{\gb}{\beta}

\newcommand{\go}{\omega}

\makeatletter
\def\captionfont@{\footnotesize}
\def\captionheadfont@{\scshape}

\long\def\@makecaption#1#2{%
  \vspace{2mm} \setbox\@tempboxa\vbox{\color@setgroup
  \advance\hsize-6pc\noindent
  \captionfont@\captionheadfont@#1\@xp\@ifnotempty\@xp
  {\@cdr#2\@nil}{.\captionfont@\upshape\enspace#2}%
  \unskip\kern-6pc\par \global\setbox\@ne\lastbox\color@endgroup}%
  \ifhbox\@ne 
  \setbox\@ne\hbox{\unhbox\@ne\unskip\unskip\unpenalty\unkern}
  \ifdim\wd\@tempboxa=\z@ 
  \setbox\@ne\hbox to\columnwidth{\hss\kern-6pc\box\@ne\hss}
  tempboxa contained more than one line
  \setbox\@ne\vbox{\unvbox\@tempboxa\parskip\z@skip
  \noindent\unhbox\@ne\advance\hsize-6pc\par}%
\fi
  \ifnum\@tempcnta<64 
    \addvspace\abovecaptionskip
    \moveright 3pc\box\@ne
  \else 
    \moveright 3pc\box\@ne
    \nobreak
    \vskip\belowcaptionskip
  \fi
\relax
}
\makeatother
\def\writefig#1 #2 #3 {\rlap{\kern #1 truecm
\raise #2 truecm \hbox{#3}}}

\begin{document}

\title[Disordered pinning models]{
Localization transition in disordered pinning models.
Effect of randomness on the critical properties.
\footnote{
Lecture Notes from the {\sl $5^{th}$ Prague Summer School on
Mathematical Statistical Mechanics}, September 11-22, 2006}
}

\author{Fabio Lucio Toninelli}
\address{
Laboratoire de Physique, UMR-CNRS 5672, ENS Lyon, 46 All\'ee d'Italie,
69364 Lyon Cedex 07, France
\hfill\break
\phantom{br.}{\it Home page:}
{\tt http://perso.ens-lyon.fr/fabio-lucio.toninelli}}
\email{fltonine@ens-lyon.fr}
\date{\today}

\begin{abstract}
These notes are devoted to the statistical mechanics of directed
polymers interacting with one-dimensional spatial defects. We are
interested in particular in the situation where frozen disorder is
present. These polymer models undergo a {\sl
localization/delocalization transition}.  There is a large
(bio)-physics literature on the subject since these systems describe,
for instance, the statistics of thermally created loops in DNA double
strands and the interaction between $(1+1)$-dimensional interfaces and
disordered walls. In these cases the transition corresponds,
respectively, to the DNA denaturation transition and to the wetting
transition.  More abstractly, one may see these models as random and
inhomogeneous perturbations of renewal processes.

The last few years have witnessed a great progress in the mathematical
understanding of the equilibrium properties of these systems.
In particular, many rigorous results about the location of the critical
point, about critical exponents and path properties of the polymer in the
two thermodynamic phases (localized and delocalized) are now available.

Here, we will focus on some aspects of this topic - in particular, on
the non-perturbative effects of
disorder.  The mathematical tools employed range from renewal theory
to large deviations and, interestingly, show tight connections with
techniques developed recently in the mathematical study of mean field
spin glasses.
\\ \\ 2000 \textit{Mathematics
Subject Classification:  60K35, 82B44, 82B41, 60K05} \\ \\
\textit{Keywords: Pinning and Wetting Models, Localization transition, Harris
Criterion, Critical Exponents, Correlation Lengths, Renewal Theory,
Interpolation and Replica Coupling}
\end{abstract}

\maketitle

\newpage
\tableofcontents

\newpage

\section{Introduction and motivations}
Consider a Markov chain $\{S_n\}_{n\in\N}$ on some state space
$\Omega$, say, $\Omega=\Z^d$. We can unfold $S$ along the discrete
time axis, i.e., we can consider the sequence $\{(n,S_n)\}_{n\in\N}$
and interpret it as the configuration of a directed polymer in the
space $\N\times \Omega$. In the examples which motivate our analysis,
the discrete time is actually better interpreted as one of the space
coordinates.  The ``directed'' character of this polymer just refers
to the fact that the first coordinate, $n$, is always increasing. In
particular, the polymer can have no self-intersections.  Some
assumptions on the law of the Markov chain will be made in Section
\ref{sec:model}, where  the model is defined precisely.  Now let $0$
be a specific point in $\Omega$, and assume that the polymer receives
a reward $\epsilon$ (or a penalty, if $\epsilon<0$) whenever $S_n=0$,
i.e., whenever it touches the {\sl defect line} $\N\times \{0\}$. In
other words, the probability of a configuration of $\{S_1,S_2,\ldots,
S_N\}$ is
modified by an exponential, Boltzmann-type factor
$$
\exp\left(\epsilon\sum_{n=1}^N {\bf 1}_{\{S_n=0\}}\right).
$$ It is clear that if $\epsilon>0$ contacts with the defect line are
enhanced with respect to the $\epsilon=0$ (or free) case, and that the
opposite is true for $\epsilon<0$. One can intuitively expect that the
in the thermodynamic limit $N\to\infty$ a phase transition occurs: for
$\epsilon>\epsilon_c$ the polymer stays close to the defect line
essentially for every $n$, while for $\epsilon<\epsilon_c$ it is
repelled by it and touches it only at a few places.  This is indeed
roughly speaking what happens, and the transition is given the name of
{\sl localization/delocalization transition}.  We warn the reader that
it is not true in general that the critical value is $\epsilon_c=0$:
if the Markov chain is transient, then $\epsilon_c>0$, i.e., a
strictly positive reward is needed to pin the polymer to the defect
line (cf. Section \ref{rem:trans}).

A more interesting situation is that where the constant
repulsion/attraction $\epsilon$ is replaced by a local,
site-dependent repulsion/attraction $\epsilon_n$. One can for instance
consider the situation where $\epsilon_n$ varies periodically in $n$,
but we will rather concentrate on the case where $\epsilon_n$ are
independent and identically distributed (IID) random variables. We
will see that, again, the transition exists when, say, the average
$\epsilon$ of $\epsilon_n$ is varied. However, in this case the
mechanism is much more subtle. This is reflected for instance in the
counter-intuitive fact that $\epsilon_c$ may be negative: a globally
repulsive defect line can attract the polymer! Presence of disorder
opens the way to a large number of exciting questions, among which we
will roughly speaking select the following one: how are the critical
point and the critical exponents influenced by disorder?

There are several reasons to study disordered pinning models:
\begin{itemize}
\item there is a vast physics and bio-physics literature on the subject,
with  intriguing (but often contradictory) theoretical predictions and
numerical/experimental observations. See also Section \ref{rem:trans};
\item they are interesting generalizations of classical renewal
sequences.  From
this point of view they raise new questions and challenges, like the
problem of the speed of convergence to equilibrium for the renewal
probability in absence of translation invariance
(cf. in particular Section \ref{sec:xi});
\item finally (and this is my main motivation) they are genuinely
quenched-disordered systems where randomness has deep,
non-perturbative effects. With respect to other systems like
disordered ferromagnets or spin glasses, moreover, disordered pinning
models have the advantage that their homogeneous counterparts are
under full mathematical control. These models, therefore, turn out to
be an ideal testing ground for theoretical physics arguments like the
Harris criterion and renormalization group analysis.
\end{itemize}
It is also quite encouraging, from the point of view of mathematical
physics, that rigorous methods have been able not only to confirm
predictions made by theoretical physicists, but in some cases also to
resolve controversies (it is the case for instance of the results in
Section \ref{sec:a>}, which disprove some claims appeared previously
in the physical literature).

\subsection{A side remark on literature and on the scope of these notes}
A excellent recent introductory work on pinning models with quenched
disorder (among other topics) is the book \cite{cf:GB} by Giambattista
Giacomin. In order to avoid the risk of producing a {\sl r\'esum\'e}
of it, we have focussed on aspects which are not (or are only
tangentially) touched in \cite{cf:GB}.  On the other hand, we will say
very little about ``polymer path properties'', to which Chapters 7 and
8 of \cite{cf:GB} are devoted. A certain degree of overlap is however
inevitable, especially in the introductory sections
\ref{sec:model} and \ref{sec:homo};
results taken from \cite{cf:GB} will be often stated without proofs (unless
they are essential in the logic of these notes).

We would also like to mention that some of the results of these notes
apply also to a model much related to disordered pinning, namely {\sl
random heteropolymers (or copolymers) at selective interfaces}. It is
the case, for instance, of the results of Sections \ref{sec:a>} and
\ref{sec:xi}. We have chosen to deal only with the pinning model for
compactness of presentation, but we invite readers interested in the
heteropolymer problem to look, for instance, at \cite{cf:BdH},
\cite{cf:Monthus}, \cite{cf:GB} and references therein.

\section{The model and its free energy}

\label{sec:model}

\subsection{The basic renewal process (``the free polymer'')}
Our starting point will be a renewal $\tau$ on the integers,
$\tau:=\{\tau_i\}_{i=0,1,2,\ldots}$, where $\tau_0=0$ and
$\{\tau_i-\tau_{i-1}\}_{i\ge1}$ are IID positive and integer-valued
random variables.  The law of the renewal will be denoted by $\bP$,
and the corresponding expectation by $\bE$. In terms of the ``directed
polymer picture'' of the introduction, $\bP$ is the law of the set
$\tau$ of the points where the polymer touches the defect line, in
absence of interaction: $\tau=\{n:S_n=0\}$ (cf. also Section
\ref{rem:trans}).  We assume that $(\tau_i-\tau_{i-1})$ or,
equivalently, $\tau_1$ is $\bP$-almost surely finite: if
\begin{eqnarray}
\label{eq:K0}
K(n):=\bP(\tau_1=n),
\end{eqnarray}
this amounts to requiring $\sum_{n\in \N} K(n)=1$. This, of course, implies
that the renewal is recurrent: $\bP$-almost surely, $\tau$ contains
infinitely many points. A second assumption is that $K(.)$ has a power-like
tail. More precisely, we require that
\begin{eqnarray}
\label{eq:K}
K(n)=\frac{L(n)}{n^{1+\alpha}}\mbox{\;\;for every\;\;}  n\in \N,
\end{eqnarray}
for some $\alpha\ge0$ and a slowly varying function $L(.)$. We recall
that a function $(0,\infty)\ni x\rightarrow L(x)\in (0,\infty)$ is
said to be slowly varying at infinity if \cite{cf:bingh}
\begin{eqnarray}
\label{eq:slow}
\lim_{x\to\infty}\frac{L(rx)}{L(x)}=1
\end{eqnarray}
for every $r>0$. In particular, a slowly varying function diverges or
vanishes at infinity slower than any power. The interested reader may
look at \cite{cf:bingh} for properties and many interesting applications of
slow variation.  Of course, every positive function $L(.)$ having a
non-zero limit at infinity is slowly varying. Less trivial examples
are $L(x)=(\log (1+x))^\gamma$ for $\gamma\in \R$.

Observe that the normalization condition $\sum_{n\in \N} K(n)=1$
implies that, if $\alpha=0$, $L(.)$ must tend to zero at infinity
(cf. also Section
\ref{rem:trans} below for an example).  

It is important to remark that typical configurations of $\tau$ are
very different according to whether $\alpha$ is larger or smaller than
$1$. Indeed the average distance between two successive points,
\begin{eqnarray}
\bE\,(\tau_i-\tau_{i-1})=\sum_{n\in \N}nK(n),
\end{eqnarray}
is finite for $\alpha>1$ and infinite for $\alpha<1$. In standard
terminology, $\tau$ is positively recurrent (i.e., $\tau$ occupies a
finite fraction of $\N$) for $\alpha>1$ and null-recurrent for
$\alpha<1$ (the density of $\tau$ in $\N$ is zero). This is a simple
consequence of the classical renewal theorem \cite[Chap. I,
Th. 2.2]{cf:asmussen}, which states that
\begin{eqnarray}
\label{eq:renewalth}
\lim_{n\to\infty}\bP(n\in\tau)=\frac1{\sum_{n\in \N}nK(n)}.
\end{eqnarray}
The distinction $\alpha\gtrless 1$ plays
an important role, especially in the behavior of the homogeneous
pinning model (cf. Section \ref{sec:homo}).  Later on we will see the
emergence of an even more important threshold value: $\alpha_c=1/2$.

\begin{rem}\rm
For $\alpha=1$, the question whether the renewal  is positively or null
recurrent
is determined by the behavior at infinity of $L(.)$: from
\eqref{eq:renewalth} we see that $\tau$ is finitely recurrent iff
$\sum_n L(n)/n<\infty$. For instance, one has null recurrence if $L(.)$ has
a positive limit at infinity.
\end{rem}

\subsection{The model in presence of interaction} Now we want to introduce
an interaction which favors the occurrence of a renewal at some points
and inhibits it at others. To this purpose,  let $\go$ (referred to
as {\sl quenched randomness} or {\sl random charges}) be a sequence
$\{\go_n\}_{n\in\N}$ of IID random variables with law $\bbP$. The
basic assumption on $\go_n$, apart from the fact of being IID, is that
$\bbE \,\go_1=0$ and $\bbE\,\go_1^2=1$.  These are rather conventions
than assumptions, since by varying the parameters $\beta$ and $h$ in
Eq. \eqref{eq:Fo} below one can effectively tune average and variance
of the charges. To be specific, in these notes we will consider only
two (important) examples: the Gaussian case $\go_1\stackrel d=
\mathcal N(0,1)$ and the bounded case, $|\go_1|\le C<\infty$. Many
results are expected (or proven) to hold in wider generality and
a few remarks in this direction are scattered throughout the notes.

We are now ready to define the free energy of our model: given $h\in\R$,
$\beta\ge0$ and $N\in\N$  let
\begin{eqnarray}
\label{eq:Fo}
F^\go_N(\beta,h):=\frac1N \log Z_{N,\go}(\beta,h):=
\frac1N \log \bE\left(e^{
\sum_{n=1}^N(\beta\go_n+h)\delta_n}\delta_N\right),
\end{eqnarray}
where for notational simplicity we put $\delta_n:=\ind_{\{n\in
\tau\}}$, ${\bf 1}_A$ being the indicator function of a set $A$. The 
 {\sl quenched average} of the free energy, or quenched free energy
 for short, is defined as
\begin{eqnarray}
\label{eq:F}
F_N(\beta,h):= \bbE F^\go_N(\beta,h).
\end{eqnarray}
Note that the factor $\delta_N$ in \eqref{eq:Fo} corresponds to
imposing the boundary condition $N\in \tau$ (the boundary condition
$0\in\tau$ at the left border is implicit in the law $\bP$). One could
equivalently work with free boundary conditions at $N$ (i.e., replace
$\delta_N$ by $1$). The infinite-volume free energy would not change,
but some technical steps in the proofs of some results would be
(slightly) more involved.

We need also a notation for the Boltzmann-Gibbs average: given
a realization $\go$ of the randomness and a system size $N$, for a
$\bP$-measurable function $f(.)$  set
\begin{eqnarray}
\label{eq:Boltz}
\bE^{\beta,h}_{N,\go}(f):=\frac{\bE\left(f(\tau)\,e^{
\sum_{n=1}^N(\beta\go_n+h)\delta_n}\delta_N\right)
}{Z_{N,\go}(\beta,h)}
\end{eqnarray}

\subsection{Existence and non-negativity of the free energy}
As usual in statistical mechanics, one is (mostly) interested in the
thermodynamic limit (i.e., the limit $N\to\infty$).  A classical
question concerns the existence of the thermodynamic limit of the free
energy, and its dependence on the realization of the randomness
$\go$. In the context of the models we are considering, the answer
is well established:
\begin{theorem} \cite[Th. 4.1]{cf:GB}
If $\bbE|\go_1|<\infty$, the limit
\begin{eqnarray}
F(\beta,h):=\lim_{N\to\infty}\frac1N\log Z_{N,\go}(\beta,h)
\end{eqnarray}
exists for every $\beta\ge0,h\in\R$ and it is $\bbP(\dd \go)$-almost
surely independent of $\go$.
\end{theorem}
Of course, the limit does depend in general on the law $\bbP$ of the
disorder.

Note that the only assumption on disorder, apart from the IID
character of the charges, is finiteness of the first moment, so that
existence and self-averaging of the infinite-volume free energy
holds in much wider generality than in the cases of Gaussian or
bounded disorder we are considering here.

Some properties of the free energy come essentially for free: in
particular, $F(\beta,h)$ is convex in $(\beta,h)$, non-decreasing in
$h$, continuous everywhere and differentiable almost everywhere as a
consequence of convexity.  Another easy fact is that the sequence
$\{N\,F_N(\beta,h)\}_{N\in\N}$ is super-additive: for every $N,M\in
N$, one has $(N+M)F_{N+M}(\beta,h)\ge N F_N(\beta,h)+M F_M(\beta,h)$.
This is easily proven:
\begin{eqnarray}
(N+M)F_{N+M}(\beta,h)&=&\bbE \log \bE\left
(e^{\sum_{n=1}^{N+M}(\beta\go_n+h)\delta_n}
\delta_{N+M}\right)\\\nonumber
&\ge&
\bbE \log \bE\left (e^{\sum_{n=1}^{N}(\beta\go_n+h)\delta_n}
\delta_N e^{\sum_{n=N+1}^{N+M}(\beta\go_n+h)\delta_n}
\delta_{N+M}\right)\\\nonumber
&=&N F_N(\beta,h)+M F_M(\beta,h),
\end{eqnarray}
where in the last step we used invariance of $\bbP$ with respect to
left shifts and the renewal property of $\bP$.  It is a standard fact
that super-additivity implies
\begin{eqnarray}
\label{eq:superadd}
F(\beta,h)\ge F_N(\beta,h)\mbox{\;\;for every\;\;}N\in\N.
\end{eqnarray}

\subsection{Contact fraction and critical point}
As we already mentioned, the interest in this class of models is mainly
due to the fact that they show a so-called {\sl localization-delocalization}
transition. This is best understood in view of the elementary bound
$F(\beta,h)\ge0$. This positivity property is immediate to prove:
\begin{eqnarray}
\label{eq:F>0}
F_N(\beta,h)\ge \frac1N\bbE \log \bE \left(e^{
\sum_{n=1}^N(\beta\go_n+h)\delta_n}\ind_{\{\tau_1=N\}}\right)=\frac hN+
\frac1N\log K(N)
\end{eqnarray}
and the claimed non-negativity in the limit follows from \eqref{eq:K}.
Recalling that $F(\beta,h)$ is non-decreasing in $h$, for a given
$\beta$ the localization/delocalization critical point is defined to
be
\begin{eqnarray}
\label{eq:hcrit}
h_c(\beta):=\sup\{h:F(\beta,h)=0\}
\end{eqnarray}
and the function $\beta\rightarrow h_c(\beta)$ is referred to as
the {\sl critical line}. The region of parameters
$$
\mathcal L:=\{(\beta,h):\beta\ge0,h>h_c(\beta)\}
$$
and
$$
\mathcal D:=\{(\beta,h):\beta\ge0,h\le h_c(\beta)\}
$$ are referred to as {\sl localized} and {\sl delocalized} phases,
respectively.  Since level sets of a convex function are convex, $\cL$
is a convex set and the function $h_c(.):[0,\infty)\ni\beta\rightarrow
h_c(\beta)$ is concave. The reason for the names ``localized'' and
``delocalized'' can be understood looking at the so-called {\sl
contact fraction} $\ell_N$, defined through
\begin{eqnarray}
\label{eq:cf}
\ell_N:=\frac{|\tau\cap \{1,\ldots,N\}|}N
\end{eqnarray}
and taking values between $0$ and $1$ (as usual, $|A|$ denotes the
cardinality of a set $A$).  It is immediate to check that
\begin{eqnarray}
\partial_h F^\go_N(\beta,h)=\bE^{\beta,h}_{N,\go}(\ell_N)
\end{eqnarray}
and, by standard arguments based on convexity, this equality survives
in the thermodynamic limit whenever the free energy is differentiable:
\begin{eqnarray}
\lim_{N\to\infty}
 \bE^{\beta,h}_{N,\go}(\ell_N)\stackrel {a.s.}=
\partial_h F(\beta,h)\;\;\;\mbox{for every\;} h\mbox{\; such that\;\;}
\partial^+_h F(\beta,h)=\partial^-_h F(\beta,h).
\end{eqnarray}
We have already mentioned that differentiability holds for
Lebesgue-almost every value of $h$. However, much more than this is
true: as it was proven in \cite{cf:GT_ALEA}, differentiability
(actually, infinite differentiability) in $h$ holds whenever
$h>h_c(\beta)$. We can therefore conclude the following: for
$h<h_c(\beta)$ (or for $h\le h_c(\beta)$ if $F(\beta,h)$ is
differentiable at $h_c(\beta)$) the thermal average of the contact
fraction tends for to zero for $N\to\infty$ (almost surely in the
disorder), while for $h>h_c(\beta)$ it tends to $\partial_h
F(\beta,h)>0$. The average contact fraction plays the role of an order
parameter, like the spontaneous magnetization in the Ising model,
which is zero above the critical temperature and positive below it.

Actually, much more refined statements about the behavior of
the contact fraction in the two phases are available. In particular:
\begin{itemize}
\item for statements concerning the localized phase we refer to
\cite{cf:GT_ALEA}. There, it is proven that, roughly speaking, not only typical
configurations $\tau$ have a number $$N\, \ell_N\sim N\,\partial_h
F(\beta,h)$$ of points, but also that these points are rather
uniformly distributed in $\{1,\ldots,N\}$: long gaps between them are
exponentially suppressed, and the largest gap is of order $\log N$
(cf. Theorem \ref{th:gap} below);
\item for $h<h_c(\beta)$ we refer to \cite{cf:GT_ptrf} and \cite[Ch. 8]{cf:GB},
where it is proven that $\ell_N$ is typically at most of order
$(\log N)/N$.
\end{itemize}
In this sense, if one goes back to the pictorial image of $\tau$ as
the set of points of polymer-defect contact, one sees that the
definition of (de)-localization in terms of free energy, as given
above, does indeed correspond to the intuitive idea in terms of path
properties: in $\cL$ the polymer stays at distance $O(1)$ from the
defect, while in $\cD$ it wanders away from it and touches it only a
small (at most $\log N$) number of times.

The reader should remark that we have made no conclusive statement
about the behavior of the contact fraction at $h_c(\beta)$, since we
have not attacked yet the very important question of the regularity of
the free energy at the critical point. This will be the subject of
Sections \ref{sec:homo} and \ref{sec:a12}.

\subsection{Quenched versus annealed free energy}
\label{sec:qa}
Inequality \eqref{eq:F>0} is a good example of how selecting a
particular subset of configurations (in that case, those for which
$\tau_1=N$) provides useful free energy lower bounds. For more refined
results in this direction we refer to \cite{cf:AS} and
\cite[Sec. 5.2]{cf:GB}.  There, this technique is employed to prove
that $h_c(\beta)$ is strictly decreasing as a function of $\beta$
which implies in particular that, since $h_c(.)$ is concave,
$h_c(\beta)$ tends to $-\infty$ for $\beta\to\infty$.  This
corresponds to the apriori non-intuitive fact that, as mentioned in
the introduction, even if the charges are on average repulsive the
defect line can pin the polymer.  This is purely an effect of spatial
inhomogeneities due to disorder: for $\beta$ large, it is convenient
for the polymer to touch the defect line in correspondence of
attractive charges, where it gets a reward $\beta \go_n+h>>1$, while
the entropic cost of avoiding the repulsive charges is independent of
$\beta$.  Free energy lower bounds were obtained also in the study of
a different model, the {\sl heteropolymer at a selective interface},
in
\cite{cf:BG}.

Free energy upper bounds are on the other hand more subtle to get. An
immediate one can be however obtained by a simple application of
Jensen's inequality:
\begin{eqnarray}
\label{eq:jensen}
F_N(\beta,h)\le \frac1N \log \bbE Z_{N,\go}(\beta,h)=\frac1N \log
\bE \left(e^{\sum_{n=1}^N(h+\log M(\beta))\delta_n}
\delta_N\right)\\\nonumber=F_N(0,h+\log M(\beta))=:F_N^a(\beta,h),
\end{eqnarray}
where $M(\beta):=\bbE\, e^{\beta \go_1}$. In particular, $\log
M(\beta)=\beta^2/2$ in the case of Gaussian disorder.  $F^a(\beta,h):=
F(0,\beta+\log M(\beta))$ is referred to as {\sl annealed free
energy}, and we see that it is just the free energy of the homogeneous
system (with the same choice of $K(.)$) computed for a shifted value
of $h$.  The physical interpretation of the annealed free energy is
clear: since configurations of $\go$ and $\tau$ are averaged on the
same footing, it corresponds to a system where impurities can
thermalize on the same time-scales as the ``polymer degrees of
freedom'' (i.e., $\tau$). This is not the physical situation one wishes
to study (quenched disorder corresponds rather to impurities which are
frozen, or which can evolve only on time-scales which are so long that
they can be considered as infinite from the experimental point of
view). All the same, the information provided by
\eqref{eq:jensen} is not at all empty.  Define first of all the {\sl
annealed critical point} as
\begin{eqnarray}
\label{eq:}
h^a_c(\beta):=\sup\{h: F^a(\beta,h)=0\}=h_c(0)-\log M(\beta).
\end{eqnarray}
Thanks to \eqref{eq:jensen} and \eqref{eq:hcrit}, one has immediately
\begin{eqnarray}
h_c(\beta)\ge h_c(0)-\log M(\beta),
\end{eqnarray}
a bound which, as will be discussed in Section \ref{sec:a<}, is optimal
for $\alpha<1/2$ and $\beta$ small.

\subsection{Back to examples and motivations}

\label{rem:trans}
Typical examples of renewal sequences satisfying \eqref{eq:K0},
\eqref{eq:K} are the following. Let $\{S_n\}_{n\ge0}$ be the simple
random walk (SRW) on $\Z$, with law $\bP^{SRW}$ and $S_0:=0$, i.e.,
$\{S_n-S_{n-1}\}_{n\in\N}$ are IID symmetric random variables with
values in $\{-1,+1\}$.  Then, it is known that \cite{cf:Feller}
$\tau:=\{n\in \N:S_{2n}=0\}$ is a null-recurrent renewal sequence such
that the law of $\tau_1$ satisfies \eqref{eq:K} with $\alpha=1/2$ and
$L(.)$ asymptotically constant. The reason why one looks only at even
values of $n$ in the definition of $\tau$ in this case is due just to
the periodicity of the SRW.  If instead one takes the SRW on $\Z^2$,
then $\tau$ (defined exactly as above) is always a null-recurrent
renewal but in this case $\alpha=0$ and $L(n)\sim c /(\log n)^2$
\cite{cf:SRW}. Note that in this case, the presence of the slowly
varying function $L(.)$ is essential in making $K(.)$ summable.

What happens in the case of the SRW on $\Z^d$ when $d\ge3$? This
example does not fall directly into the class we are considering since
this process is transient, and therefore the set $\tau$ of its returns
to zero is a transient renewal sequence. However this is not too
bad. Indeed, suppose more generally that one is given $K(.)$ which
satisfies
\eqref{eq:K} but such that $\Sigma:=\sum_{n\in\N} K(n)<1$, i.e.,
$K(.)$ is a sub-probability on $\N$. Then, one may define $\hat
K(n):=K(n)/\Sigma$ which is obviously a probability. It is easy to
realize from Eq. \eqref{eq:Boltz} that the Gibbs measure (and free
energy) of the model defined starting from $K(.)$ is the same as that
obtained starting from $\hat K(.)$, provided that $h$ is replaced by
$h+\log \Sigma$.  The case where $\tau$ are the zeros of the SRW on
$\Z^d$ with $d\ge3$ can then be included in our discussion:
Eq. \eqref{eq:K} holds with $\alpha=d/2-1$ and $L(.)$ asymptotically
constant.  In the following we will therefore always assume, without
loss of generality, that $\tau$ is recurrent.

We conclude this section by listing a couple of examples of
(bio)-physical situations where disordered pinning models are
relevant:
\begin{itemize}
\item {\sl Wetting of $(1+1)$-dimensional disordered substrates}
\cite{cf:Derrida} \cite{cf:Fogacs}.  Consider a two-dimensional system
at a first order phase transition, e.g., the $2d$-Ising model at zero
magnetic field and $T<T_c$, or a liquid-gas system on the coexistence
line.  Assume that the system is enclosed in a square box with
boundary conditions imposing one of the two phases along the bottom
side of the box and the other phase along the other three sides.  For
instance, for the Ising model one can impose $+$ boundary conditions
(b.c.)  along the bottom side and $-$ b.c. along the other ones; for
the liquid-gas model, one imposes that the bottom of the box is in
contact with liquid and that side and top walls are in contact with
gas. Then, there is necessarily an interface joining the two bottom
corners of the box and separating the two phases.  At very low
temperature, it is customary to describe this interface as a
one-dimensional symmetric random walk (not necessarily the SRW)
conditioned to be non-negative, the non-negativity constraint
reflecting the fact that the interface cannot exit the box. The
directed character of the random walk implies in particular that one
is neglecting the occurrence of bubbles or overhangs in the
interface. An interesting situation occurs when the bottom wall is
``dirty'' and at each point has a random interaction with the
interface: at some points the wall prefers to be in contact with the
gas (or $-$ phase), and therefore tries to pin the interface, while at
other points it prefers contact with the liquid (or $+$ phase) and
repels the interface. Of course, this non-homogeneous interaction is
encoded in the charges $\go_n$.  In this context, the
(de)-localization transition is called {\sl wetting transition}. This
denomination is clear if we think of the liquid-gas model: the
localized phase corresponds to an interface which remains at finite
distance from the wall (the wall is dry), while in the delocalized
phase there are few interface-wall contacts and the height of the
liquid layer on the wall diverges in the thermodynamic limit: the wall
is wet.  It is known that, in great generality \cite{cf:Feller}, the
law of the first return to zero of a one-dimensional random walk
conditioned to be non-negative is of the form \eqref{eq:K} with
$\alpha=1/2$ and $L(.)$ asymptotically constant (this process is
transient but this fact is not so relevant, in view of the discussion
at the beginning of the present section).
\item {\sl Formation of loops under thermal excitation and
denaturation of DNA molecules} in the Poland-Scheraga (PS)
approximation \cite{cf:Cule}. Neglecting its helical structure, the
DNA molecule is essentially a double strand of complementary units,
called ``bases''.  Upon heating, the bonds which keep base pairs
together can break and the two strands can partly or entirely separate
(cf. figure below). This separation, or denaturation, can be described
in the context of our disordered pinning models. The set $\tau$
represents the set of bases whose bond is not broken.  In the
localized phase $\tau$ contains $O(N)$ points ($N$ being interpreted
here as the total DNA length), i.e., corresponds to the phase where
the two strands are still essentially tightly bound. In the
delocalized (or denaturated) phase, on the contrary, only few bases
pairs are bound.  In formulating the PS model, one usually takes a
value $\alpha\simeq 2.12$ (cf. \cite{cf:kafri} for a justification of
this choice) and (in our notations, which are not necessarily those of
the literature on the PS model)
$$
L(n)=\sigma \mbox{\;\;for \;\;} n\ge 2,
$$ where $\sigma$ (the {\sl cooperativity parameter}) is a small
number, usually of the order $10^{-5}$, while $L(1)$ is fixed by the
normalization condition $\sum_{n\in\N}K(n)=1$.  Quenched disorder
corresponds here to the fact that bases of the different types are
placed inhomogeneously along the DNA chain.  We refer to \cite[Section
1.4]{cf:GB} for a very clear introduction to the denaturation problem
and the Poland-Scheraga model.  Here we wish to emphasize only that
the renewal process $\tau$ described by such a $K(.)$ is not in
general the set of returns of a Markov chain, as it happens for instance
in the case of the wetting model described above.
\end{itemize}
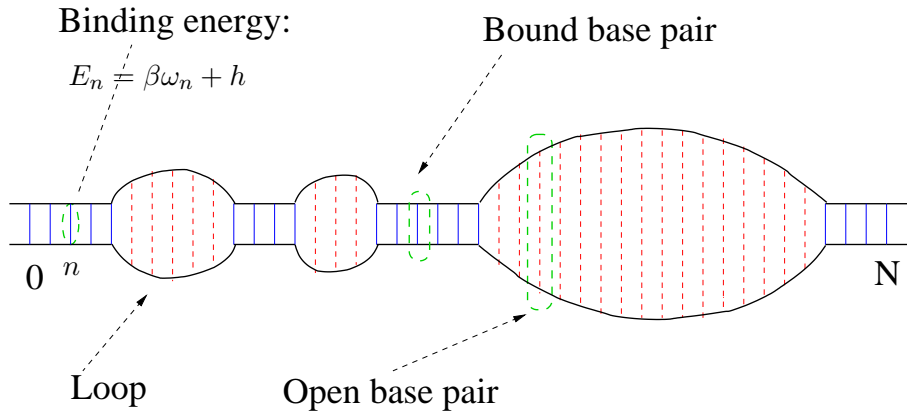
\begin{figure}[h]
\begin{center}~
\leavevmode
\epsfxsize =12 cm
\psfragscanon
\psfrag{EE}[c]{{ $E_n=\beta\go_n+h$}}
 \psfrag{n}[c]{{ $n$}}
\epsfbox{DNA.eps}
\end{center}
\end{figure}

\section{The questions we are interested in}

The main questions which will be considered in these notes are the
following:
\begin{enumerate}
\item When is the annealed bound \eqref{eq:jensen} a good one, i.e., when are
quenched and annealed systems similar?  We will see that quenched and
annealed free energies never coincide, except in the (trivial) case
where the annealed free energy is zero (i.e., the annealed model is
delocalized). However, this does not mean that the solution of the
annealed system gives no information about the quenched one. For
instance we will show that, for $\alpha<1/2$ and weak enough disorder,
the quenched critical point coincides with the annealed one.  This
will be discussed in Section \ref{sec:a<}.
\item What is the order of the transition? Critical exponents
(in particular, the {\sl specific heat exponent}, cf. next section)
can be exactly computed for the homogeneous model. The Harris criterion
predicts that for small $\beta$ critical exponents are those of the
$\beta=0$ (or annealed) model if $\alpha<1/2$, and are different if
$\alpha>1/2$.  This is the question of {\sl disorder relevance},
discussed in Sections
\ref{sec:a<}-\ref{sec:a>}.
\item Truncated correlations functions are known to decay exponentially
at large distance, in the localized phase. What is the behavior of the
correlation length when the transition is approached? We will see
that, due to the presence of quenched disorder, one can actually
define two different correlation lengths. In specific cases, we will
identify these correlation lengths and give bounds on the critical
exponents which govern their divergence at $h_c(\beta)$.
\end{enumerate}

\section{The homogeneous model}
\label{sec:homo}
In absence of disorder ($\beta=0$)
the model is under full mathematical control; in particular,
critical point and the order of the transition can be computed exactly.
In this section, we collect a number of known results, referring to
\cite[Chapter 2]{cf:GB} for their proofs.

The basic point is that the free energy $F(0,h)$ is determined as follows
\cite[Appendix A]{cf:GT_smooth}:
if the equation
\begin{eqnarray}
\label{eq:id3}
  \sum_{n\in\N}e^{-b n}K(n)=e^{-h}
\end{eqnarray}
has a positive solution $b=b(h)>0$ then $F(0,h)=b(h)$. Otherwise,
$F(0,h)=0$.  From this (recall the normalization condition
$\sum_{n\in\N}K(n)=1$), one finds immediately that $h_c(0)=0$. The
behavior of the free energy in the neighborhood of $h_c(0)$ can be
also obtained from \eqref{eq:id3}. Care has to be taken since a naive
expansion of left- and right-hand sides of \eqref{eq:id3} for $b$ and
$h$ small does not work in general. However, this analysis can be
performed without much difficulty and one can prove the following:
\begin{theorem} \cite[Th. 2.1]{cf:GB}
\label{th:hom}
  \begin{enumerate}
    \item If $\alpha=0$, $F(0,h)$ vanishes faster than any power of
$h$ for $h\searrow0$.
\item If $0<\alpha<1$ then  for $h>0$
\begin{eqnarray}
\label{eq:asint}
  F(0,h)=h^{1/\alpha} \hat L(1/h),
\end{eqnarray}
where $\hat L(.)$ is the slowly varying function
\begin{eqnarray}
\hat L(1/h)=\left(\frac{\alpha}{\Gamma(1-\alpha)}\right)^{1/\alpha}
h^{-1/\alpha}R_\alpha(h)
\end{eqnarray}
and $R_\alpha(.)$ is asymptotically equivalent to the inverse of the
map $x\rightarrow x^\alpha L(1/x)$.
\item If $\alpha=1$ and $\sum_{n\in\N}n K(n)=\infty$ then
$F(0,h)=h\,\hat L(1/h)$ for some slowly varying function
$\hat L(.)$ which vanishes at infinity.
\item  If $\sum_{n\in\N}n K(n)<\infty$ (in particular, if $\alpha>1$)
  \begin{eqnarray}
    F(0,h)\stackrel{h\searrow0}{\sim} \frac h{\sum_{n\in\N}n K(n)}.
  \end{eqnarray}
  \end{enumerate}
\end{theorem}

In particular, note that in the situation (4), i.e., if $\tau$ is
positively recurrent under $\bP$, the transition is of first order:
the free energy is not differentiable at $h_c(0)=0$, i.e., the average
contact fraction has a finite jump in the thermodynamic limit.  This
is analogous to what happens for the Ising model in dimension $d\ge2$:
if $T<T_c$ and one varies the magnetic field $H$ from $0^-$ to $0^+$, the
spontaneous magnetization has a positive jump and the free energy is
not differentiable.  The transition is, on the other hand, continuous
({\sl at least of second order}) if $\bP$ is the law of a
null-recurrent renewal $\tau$ and it becomes smoother as $\alpha$
decreases. In thermodynamical language, one can say that the
delocalization transition is of $k^{th}$ order ($F(\beta,.)$ is of
class $C^{k-1}$ but not of class $C^k$) for $\alpha\in (1/k, 1/{k-1})$
and of infinite order for $\alpha=0$.\footnote{ In order to decide
between $k^{th}$ and $(k+1)^{th}$ order for $\alpha=1/k$ one needs to
look also at the slowly varying function $L(.)$, as is already clear
from points (3) and (4) in the case of $k=1$. In any case, the precise
statement is that of Theorem \ref{th:hom}.}

In the physics literature one introduces usually the {\sl specific heat
critical exponent} $\nu$ as\footnote{\label{nota:ocio} the symbol $\nu$
for the specific heat exponent is not standard in the literature,
but we have already used the letter $\alpha$ for another purpose. The
same remark applies to the symbols we use for other critical exponents.}
\begin{eqnarray}
\label{eq:nu}
\nu=2-\lim_{h\searrow h_c(\beta)} \frac{\log F(\beta,h)}{\log (h-h_c(\beta))}
\end{eqnarray}
(provided the limit exists) and of course $\nu$ can depend on
$\beta$. From Theorem \ref{th:hom} we see that, in absence of
disorder,
\begin{eqnarray}
\label{eq:nb0}
\nu(\beta=0)=2-\max(1,1/\alpha).
\end{eqnarray} In particular, note that $\nu(\beta=0)>0$
as soon as $\alpha>1/2$ (this observation will become interesting in
the light of the results of Section \ref{sec:a>}).

 \section{Relevance or irrelevance of disorder?}
\label{sec:a12}
We have just seen that the phase transition of the homogeneous pinning
model can be of any given order - from first to infinite - depending
on the choice of $K(.)$ in \eqref{eq:K} and, in particular, on the
value of $\alpha$. In this section we discuss the effect of
disorder on the transition and we are primarily interested in the
question of disorder relevance. There are actually two distinct (but
inter-related) aspects in this question:
\begin{enumerate}
\item[{\bf Q1}]
does an arbitrarily small quantity of disorder change the critical
exponent $\nu$ (i.e., the order of the transition)?
\item[{\bf Q2}]  does the quenched critical point differ from the annealed one
for very weak disorder?
\end{enumerate}
One expects the answer to both questions to be ``no'' if $\alpha<1/2$
and ``yes'' if $\alpha>1/2$, while the case $\alpha=\alpha_c=1/2$ is more
subtle and not clear even heuristically \cite{cf:Derrida, cf:Fogacs}
(see, however, Theorem \ref{th:a32}).

The plan is the following: we will first of all (Section
\ref{sec:Harris}) make a non-rigorous computation, in the spirit of the
Harris approach \cite{cf:Harris}, which shows why the watershed value
for $\alpha$, distinguishing between relevance and irrelevance, is
expected to be $ \alpha_c=1/2$, i.e., the value for which the critical
exponent $\nu$ vanishes for the homogeneous model
(cf. \eqref{eq:nb0}).  Next, in Section
\ref{sec:UBq} we prove an upper bound for the free energy which
strictly improves the annealed bound \eqref{eq:jensen}. In the proof
of this bound we introduce the technique of {\sl interpolation}, by
now classical in spin glass theory but sort of new in this context. We
would like to emphasize that interpolation (and replica coupling,
cf. Section \ref{sec:pfs}) techniques have proven recently to be
extremely powerful in the analysis of mean field spin glass models,
cf. for instance \cite{cf:GuT_cmp}, \cite{cf:Ass}, \cite{cf:T}, while
their relevance in the domain of disordered pinning model had not
been realized  clearly so far.

As a byproduct, our new upper bound partially justifies the heuristic
expansion of Section
\ref{sec:Harris}. The question of relevance is taken up more seriously
in Sections
\ref{sec:a<} to \ref{sec:a>}. In the former we will see, among other
results, that answers to both {\bf Q1} and {\bf Q2} are actually
``no'' for $\alpha<\alpha_c$.  In the latter, on the other hand, we
show that critical exponents {\sl are} modified by disorder for
$\alpha>\alpha_c$: in particular, we will see that
$\nu\le0$ whenever $\beta>0$.

In the whole of Section \ref{sec:a12} we consider only the case of
Gaussian disorder. This allows for technically simpler proofs, but
results can be generalized for instance to the bounded disorder case.

\subsection{Harris criterion and the emergence of $\alpha_c=1/2$}
\label{sec:Harris}
Let us note for clarity that, putting together the discussion of
Section \ref{sec:homo} and Eq. \eqref{eq:}, in the Gaussian case the
annealed critical point equals $h_c^a(\beta)=-\beta^2/2$.  The first
step of our heuristic argument is rigorous and, actually, an immediate
identity:
\begin{eqnarray}
\label{eq:id1}
F_N(\beta,h)=F_N^a(\beta,h)+\frac1N\bbE \log
\med{e^{\sum_{n=1}^N(\beta\go_n-\beta^2/2)\delta_n}}_{N,h-h_c^a(\beta)},
\end{eqnarray}
where $\med{.}_{N,h}:=\bE^{0,h}_{N,0}(.)$ is just the Boltzmann
average for the homogeneous system (cf. Eq. \eqref{eq:Boltz}).
Identity \eqref{eq:id1} can be rewritten in a more suggestive way if
we recall the last equality in \eqref{eq:jensen} and we let
$h=h^a_c(\beta)+\Delta$ with $\Delta\ge0$:
\begin{eqnarray}
\label{eq:id2}
F_N(\beta,h^a_c(\beta)+\Delta)=F_N(0,\Delta)+
R_{N,\Delta}(\beta):=F_N(0,\Delta)+\frac1N\bbE \log
\med{e^{\sum_{n=1}^N(\beta\go_n-\beta^2/2)\delta_n}}_{N,\Delta}.
\end{eqnarray}
Irrelevance of disorder amounts to the fact that, for $\beta$
sufficiently small, the ``error term'' $R_{N,\Delta}(\beta)$ is
negligible with respect to the ``main term'' $F_N(0,\Delta)$. As we
will see, the question is subtle since we are interested in both
$\Delta$ and $\beta$ small, and the two limits do not in general
commute. For the moment, let us
proceed without worrying about rigor and let us expand naively
$R_{N,\Delta}(\beta)$ for $\beta$ small and $\Delta,N$ fixed:
\begin{eqnarray}
\med{e^{\sum_{n=1}^N(\beta\go_n-\beta^2/2)\delta_n}}_{N,\Delta}=
1+\sum_{n=1}^N(\beta\go_n-\beta^2/2)\med{\delta_n}_{N,\Delta}+
\frac{\beta^2}2\sum_{n,m=1}^N\go_n\go_m\med{\delta_n\delta_m}_{N,\Delta}+
O(\beta^3).
\end{eqnarray}
Expanding the logarithm and using the fact that $\bbE\,\go_n=0$ and $\bbE
(\go_n\go_m)={\bf 1}_{\{n=m\}}$ one has, always formally,
\begin{eqnarray}
R_{N,\Delta}(\beta)=-\frac{\beta^2}{2N}\sum_{n=1}^N\left(
\med{\delta_n}_{N,\Delta}\right)^2+O(\beta^3).
\end{eqnarray}
In the limit $N\to\infty$ one has by definition of the homogeneous model
$$
\lim_{N\to\infty}\med{\ell_N}_{N,\Delta}=\lim_{N\to\infty}
\frac1N \sum_{n=1}^N\med{\delta_n}_{N,\Delta}=\partial_\Delta F(0,\Delta).
$$
Since $\med{\delta_n}_{N,\Delta}$ should not depend on $n$ as soon
as $1\ll n\ll N$, one can expect (actually, this can be proven without
much difficulty) that
\begin{eqnarray}
\limtwo{N\to\infty}{n/N\to m\in(0,1)}\med{\delta_n}_{N,\Delta}=\partial_\Delta F(0,\Delta).
\end{eqnarray}
In conclusion, we find
\begin{eqnarray}
\label{eq:expa}
F(\beta,h_c^a(\beta)+\Delta)=F(0,\Delta)-\frac{\beta^2}2
(\partial_\Delta F(0,\Delta))^2+ O(\beta^3).
\end{eqnarray}
Even without trying (for the moment) to justify this expansion or to
look more closely at the $\Delta$-dependence of the error term
$O(\beta^3)$, we can extract something important from
Eq. \eqref{eq:expa} . We know from Theorem \ref{th:hom} that, for
$\alpha<1$ and $\Delta>0$ small, $F(0,\Delta)\simeq
\Delta^{1/\alpha}$ which implies (cf. the proof of Eq. \eqref{eq:L()}
for details) that $\partial_\Delta F(0,\Delta)
\simeq \Delta^{1/\alpha-1}$. Then we see immediately that, indeed, for
$\alpha<1/2$
\begin{eqnarray}
\label{eq:lhsrhs}
\frac{\beta^2}2
(\partial_\Delta F(0,\Delta))^2\ll F(0,\Delta)
\end{eqnarray}
if $\Delta$ and $ \beta$ are small. In terms of the Harris criterion,
disordered is said to be irrelevant in this case and one can hope that
the expansion can be actually carried on at higher orders.  For $1/2<
\alpha<1$, however, this is false: even if $\beta$ is small, choosing
$\Delta$ sufficiently close to zero the left-hand side of
\eqref{eq:lhsrhs} is much larger than the right-hand side. This means
that ``disorder is relevant'' and the small-disorder expansion breaks
down immediately. The same holds for $\alpha\ge 1$, when
$F(0,\Delta)\simeq \Delta$ and $\partial_\Delta F(0,\Delta)\sim
const$. The threshold value $\alpha_c=1/2$ is clearly a ``marginal
case'' where relevance or irrelevance of disorder cannot be decided
(even on heuristic grounds) by a naive expansion in $\beta$.

The rest of this section will be devoted to give rigorous bases to
this suggestive picture. As a byproduct we will learn something
interesting for the case $1/2<\alpha<1$: while disorder is relevant
and changes the exponent $\nu$, it modifies the transition
only ``very close'' to the critical point (cf. Theorem \ref{th:a322}).

\subsection{A rigorous approach:  interpolation and an
improvement upon annealing}
\label{sec:UBq}
In Section \ref{sec:qa} we saw that a simple application of Jensen's
inequality implies $F(\beta,h)\le F^a(\beta,h)$. Here we wish to show
that this inequality is strict as soon as disorder is present
($\beta>0$) and the annealed system is localized.  Moreover, we will
partly justify the small-$\beta$ expansion of Section \ref{sec:Harris}  for
$\alpha<1/2$, showing that it provides an upper bound for the quenched
free energy.

More precisely:
\begin{theorem}\cite[Th. 2.6]{cf:T_qrc}
\label{th:rs}
For every $\beta>0$, $\alpha\ge0$ and $\Delta>0$
\begin{eqnarray}
\label{eq:RS}
F(\beta,h^a_c(\beta)+\Delta)\le \inf_{0\le q\le \Delta/\beta^2}
\left(\frac{\beta^2q^2}2+F(0,\Delta-\beta^2q)\right)<
F(0,\Delta)=F^a(\beta,h).
\end{eqnarray}
In particular, if $0\le\alpha<1/2$ there exist
constants $\beta_0>0, \Delta_0>0$
such that
\begin{eqnarray}
\label{eq:svil}
F(\beta,h^a_c(\beta)+\Delta)\le F(0,\Delta)-\frac{\beta^2}2
\left(\partial_\Delta F(0,\Delta)\right)^2(1+O(\beta^2))
\end{eqnarray}
for $\beta\le \beta_0, \Delta\le \Delta_0$, where $O(\beta^2)$ is does
not depend on $\Delta$.  On the other hand, if $\beta=0$ or $\Delta\le
0$, then $F(\beta,h_c^a(\beta)+\Delta)=F^a(\beta,
h_c^a(\beta)+\Delta)$.
\end{theorem}
About the possibility of pushing the upper bound \eqref{eq:svil} to
order higher than $\beta^2$ see Remark 3.1 in \cite{cf:T_qrc}.
It
is obvious that \eqref{eq:svil} cannot hold for $\alpha>1/2$ since, as
already observed after Eq. \eqref{eq:lhsrhs}, the right-hand side is
negative for $\Delta$ sufficiently small.

Readers familiar with mean field spin glass models will remark a
certain similarity between the variational bound \eqref{eq:RS} and the
``replica symmetric'' variational bound \cite{cf:Gu_rs} for the free
energy of the Sherrington-Kirkpatrick model. However, we do not see
a natural way to generalize \eqref{eq:RS} to include ``replica symmetry
breaking'' in analogy with \cite{cf:Gu} \cite{cf:Ass}.

{\sl Proof of Theorem \ref{th:rs}}. The proof is rather instructive
because it allows us to introduce the technique of ``interpolation'',
which will play a major role in the next subsection.  We start from
identity \eqref{eq:id2} and, for $\Delta>0, q\in \R $ and $0\le
t\le1$, we define
\begin{eqnarray}
\label{eq:Rt}
  R_{N,\Delta}(t,\beta,q):=\frac1N\bbE\log
\med{e^{\sum_{n=1}^N[\beta\sqrt t\go_n
-t\beta^2/2+\beta^2q(t-1)]\delta_n}}_{\Delta,N}.
\end{eqnarray}
In spin glass language, this would be called an ``interpolating free energy'',
since by varying the parameter $t$ it relates in a smooth way the quantity
we wish to estimate at $t=1$,
\begin{eqnarray}
\label{eq:bc1}
  R_{N,\Delta}(t=1,\beta,q)=R_{N,\Delta}(\beta)
\end{eqnarray}
to something easy at $t=0$:
\begin{eqnarray}
\label{eq:bc0}
  R_{N,\Delta}(t=0,\beta,q)=F_N(0,\Delta-\beta^2q)-F_N(0,\Delta).
\end{eqnarray}
A priori, there is no reason why $R_{N,\Delta}(t,\beta,q)$ should be
any easier to compute for $0<t<1$ than for $t=1$.  What helps us is
that the $t$-derivative of $R_{N,\Delta}(t,\beta,q)$ can be bounded
above by throwing away a (complicated) term which, luckily, has a
negative sign. To see this we need first of all manageable notations
and we will set
\begin{eqnarray}
\label{eq:manag}
\med{g(\tau)}_{N,\Delta,t}:=\frac{\med
{g(\tau) e^{\sum_{n=1}^N[\beta\sqrt t\go_n
-t\beta^2/2+\beta^2q(t-1)]\delta_n}}_{\Delta,N}}
{\med{e^{\sum_{n=1}^N[\beta\sqrt t\go_n
-t\beta^2/2+\beta^2q(t-1)]\delta_n}}_{N,\Delta}}
\end{eqnarray}
for every measurable function $g(\tau)$.
We find then
\begin{eqnarray}
\label{eq:derivq}
\frac{\dd}{\dd t}R_{N,\Delta}(t,\beta,q)=
\frac{\beta^2}N \left(-\frac12+q\right)\sum_{m=1}^N \bbE
\med{\delta_m}_{N,\Delta,t}
+\frac{\beta}{2\sqrt t N}\sum_{m=1}^N\bbE\, \go_m\med{\delta_m}_{N,\Delta,t}.
\end{eqnarray}
The last term of \eqref{eq:derivq} can be rewritten using the
Gaussian integration by parts formula
\begin{eqnarray}
\label{eq:byparts}
  \bbE \left(\go f(\go)\right)=\bbE f'(\go),
\end{eqnarray}
which holds (if $\go$ is a Gaussian random variable $\mathcal N(0,1)$) for
every differentiable function $f(.)$ such that
$\lim_{|x|\to\infty}\exp(-x^2/2)f(x)=0$. In our case,
the function $f$ is of course $\med{\delta_m}_{N,\Delta,t}$ and one finds
\begin{eqnarray}
\label{eq:derivq2}
\frac{\beta}{2\sqrt t N}\sum_{m=1}^N\bbE\, \go_m\med{\delta_m}_{N,\Delta,t}=
\frac{\beta^2}{2N}\sum_{m=1}^N \bbE \left(\med{\delta_m}_{N,\Delta,t} -
 \left(\med{\delta_m}_{N,\Delta,t}
\right)^2\right).
\end{eqnarray}
The positive term comes from the differentiation of the numerator of
$\med{\delta_m}_{N,\Delta,t}$ (recall the definition \eqref{eq:manag})
 and the negative one from the denominator, and we used the obvious
$\delta_m=(\delta_m)^2$. Putting together Eqs. \eqref{eq:derivq} and
\eqref{eq:derivq2} one has therefore
\begin{eqnarray}
\label{eq:derivq3}
\frac{\dd}{\dd t}R_{N,\Delta}(t,\beta,q)=\frac{\beta^2q^2}2-
\frac{\beta^2}{2N}\sum_{n=1}^N\bbE \left\{\left(\med{\delta_n}_{N,\Delta,t}
-q\right)^2\right\}\le \frac{\beta^2q^2}2.
\end{eqnarray}
At this point we are done: we integrate on $t$ between $0$ and $1$ inequality
\eqref{eq:derivq3}, we recall the boundary conditions \eqref{eq:bc0} and
\eqref{eq:bc1} and we get
\begin{eqnarray}
R_{N,\Delta}(\beta)\le
F_N(0,\Delta-\beta^2q)-F_N(0,\Delta)+\frac{\beta^2q^2}2.
\end{eqnarray}
Together with Eq. \eqref{eq:id1}, taking $N\to\infty$ limit and
minimizing over $q$ proves \eqref{eq:RS}.  Let us remark that
minimizing over $q\in \R$ or on $0\le q\le \Delta/\beta^2$ is clearly
equivalent. The strict inequality in \eqref{eq:RS} is just due to the
fact that the derivative with respect to $q$ of the quantity to be
minimized, computed at $q=0$, is negative.

The expansion \eqref{eq:svil} is just a consequence of
\eqref{eq:RS}. Remark first of all that, at the lowest order in $\beta$,
the minimizer in \eqref{eq:RS} is $q=q_\Delta:=\partial_\Delta
F(0,\Delta)$. Then, from
identity \eqref{eq:id3} one finds that there exist slowly varying functions
$L^{(i)}(.), i=1,2$ such that for $\alpha<1/2$ and $\Delta>0$
\begin{eqnarray}
\label{eq:L()}
\partial_\Delta F(0,\Delta)=\Delta^{(1-\alpha)/\alpha} L^{(1)}(1/\Delta),
\;\;\;\;
\partial^2_\Delta F(0,\Delta)=\Delta^{(1-2\alpha)/\alpha}L^{(2)}(1/\Delta).
\end{eqnarray}
Let us show for instance the first equality. Differentiating both sides of
\eqref{eq:id3} with respect to $\Delta$ one finds
\begin{eqnarray}
\partial_\Delta F(0,\Delta)=\frac{e^{-\Delta}}{\sum_{n\in\N}
n^{-\alpha}L(n)\exp(-F(0,\Delta)n)}.
\end{eqnarray}
Using Theorems \ref{th:A1} and \ref{th:A2} one has then, for $\Delta\to0$
(i.e., for $F(0,\Delta)\to0$)
\begin{eqnarray}
\partial_\Delta F(0,\Delta)\stackrel{\Delta\searrow0}\sim \frac{
\Gamma(2-\alpha)L(1/F(0,\Delta))}{(1-\alpha)F(0,\Delta)^{1-\alpha}}
\end{eqnarray}
which, together with \eqref{eq:asint}, proves the first equality in
\eqref{eq:L()} for a suitable $L^{(1)}(.)$. Note, by the way, that
thanks to \eqref{eq:L()} one has $q_\Delta<\Delta/\beta^2$ for
$\Delta,\beta$ sufficiently small (and $\alpha<1/2$, of
course). Another consequence of \eqref{eq:L()} is that
$\partial^2_\Delta F(0,\Delta)$ is bounded above by a finite constant
$C$ for, say, $\Delta\le1$. Then, a Taylor expansion gives
$$
F(0,\Delta-\beta^2 q_\Delta)\le F(0,\Delta)-\beta^2(\partial_\Delta
F(0,\Delta))^2+C\beta^4 (\partial_\Delta F(0,\Delta))^2,
$$
whence Eq. \eqref{eq:svil}.

Finally, the statement for $\beta=0$ or $\Delta\le0$ is trivial: for
$\beta=0$ there is no disorder to distinguish between quenched an
annealed free energies, and for $\Delta\le0$ one has
$F^a(\beta,h_c^a(\beta)+\Delta)=0$ which, together with
\eqref{eq:jensen} and $F(\beta,h)\ge0$, implies the statement.

\qed

\subsection{Irrelevance of disorder for $\alpha<1/2$ via replica coupling}
\label{sec:a<}
We want to say first of all that, if $0<\alpha<1/2$ and $\beta$ is
sufficiently small (i.e., if disorder is sufficiently weak),
then $h_c(\beta)=h^a_c(\beta)$.  Recalling that
$F^a(\beta,h_c^a(\beta)+\Delta)=F(0,\Delta)$, this follows immediately from
\begin{theorem}
\cite{cf:Ken, cf:T_qrc}
\label{th:alpha<32}
  Assume that either $0<\alpha<1/2$ or that
  \begin{eqnarray}
 \label{eq:orthat}
  \alpha=1/2 \mbox{\;\;and\;\;}
\sum_{n\in \N}n^{-1}L(n)^{-2}<\infty.
  \end{eqnarray}
Then, for every
$\epsilon>0$ there exist $\beta_0(\epsilon)>0$ and
$\Delta_0(\epsilon)>0$ such that, for every $\beta\le
\beta_0(\epsilon)$ and $0<\Delta<\Delta_0(\epsilon)$, one has
\begin{eqnarray}
\label{eq:risulta1}
(1-\epsilon)F(0,\Delta)\, \le \, F(\beta,h_c^a(\beta)+\Delta)\le
F(0,\Delta).
\end{eqnarray}
\end{theorem}
\medskip
Observe that this implies in particular that, under the assumptions of
the theorem, the exponent $\nu$ equals $2-1/\alpha$ as in the
homogeneous case. Indeed note that, for $\Delta$ small,
\begin{eqnarray}
\frac{\log (1-\epsilon)+\log F(0,\Delta)}{\log \Delta}\ge
\frac{\log F(\beta,h_c(\beta)+\Delta)}{\log \Delta}\ge
\frac{\log F(0,\Delta)}{\log \Delta}
\end{eqnarray}
and the statement follows taking the limit
$\Delta\to0$ from definition \eqref{eq:nu} of the
specific heat exponent.

We will see in Section \ref{sec:a>} that the same cannot hold for
$\alpha>1/2$: in that case, $\nu$ is necessarily non-positive in for
the quenched system presence of disorder, while it is positive for the
annealed system. One could therefore think that quenched and annealed
behaviors are completely different. This is however not completely
true. Indeed, the next theorem shows that $F(\beta,h)$ and
$F^a(\beta,h)$ are very close, provided that $1/2\le \alpha<1$ if one
is not too close to the critical point. More precisely one has
\begin{theorem}
\label{th:a322}
   Assume that $1/2<\alpha<1$. There exists a slowly
varying function $ \check L(.)$ and, for every $\epsilon>0$, constants
$a_1(\epsilon)< \infty$ and $\Delta_0(\epsilon)>0$ such that, if
\begin{eqnarray}
\label{eq:condiz}
a_1(\epsilon)\beta^{2\alpha/(2\alpha-1)}\check L (1/\beta)
\, \le\,  \Delta \,
\,
\le\,  \Delta_0(\epsilon),
\end{eqnarray}
the inequalities \eqref{eq:risulta1} hold.
\end{theorem}
\medskip

To see more clearly what this says on the relation
between quenched and annealed critical points, forget about the slowly
varying functions; then,  Theorem \ref{th:a322} implies
$$
0\le h_c(\beta)-h_c^a(\beta)\lesssim \beta^{2\alpha/(2\alpha-1)}.
$$
Since $2\alpha/(2\alpha-1)>2$, this
shows in particular that
\begin{eqnarray}
\lim_{\beta\searrow 0}\frac{h_c(\beta)}{h^a_c(\beta)}=1.
\end{eqnarray}

\begin{rem}\rm
\label{r:54}
Theorem \ref{th:a322} was proven in \cite[Th. 3]{cf:Ken} and then in
\cite[Th. 2.2]{cf:T_qrc}. The two results differ only  in the form of the
slowly varying function $\check L(.)$. In general,
the function  $\check L(.)$ which pops out from the proof in
\cite[Th. 2.2]{cf:T_qrc} is larger (i.e., worse) than that of
\cite[Th. 3]{cf:Ken}.
\end{rem}

Finally, we consider the ``marginal case'' $\alpha=\alpha_c=1/2$ and
$\sum_n (L(n))^{-2}n^{-1}=\infty$. This is the case, for instance, if
$\bP$ is the law of the returns of a one-dimensional symmetric random
walk, where $L(.)$ is asymptotically constant, as mentioned in Section
\ref{rem:trans}. As we mentioned, this case is still debated even in
the physical literature. The ``most likely'' scenario
\cite{cf:Derrida} is that disorder is ``marginally relevant'' in this
case: $h_c(\beta)\ne h^a_c(\beta)$ for every positive $\beta$, but the
two critical points are equal at every order in a weak-disorder
perturbation theory. Other works, e.g. \cite{cf:Fogacs}, claim on the
other hand that disorder is irrelevant in this situation.

What one can prove for the moment is the following:
\begin{theorem}\cite{cf:Ken,cf:T_qrc}
\label{th:a32}
   Assume that $\alpha=1/2$ and $\sum_{n\in \N}n^{-1}L(n)^{-2}=\infty$.
Let $\ell(.)$ be the slowly varying function (diverging at infinity)
defined by
\begin{eqnarray}
\label{eq:lN}
\sum_{n=1}^N \frac 1{n L(n)^2}\stackrel{N\to\infty}{\sim}\ell(N).
\end{eqnarray}
For every $\epsilon>0$ there exist constants
$a_2(\epsilon)<\infty$ and $\Delta_0(\epsilon)>0$
such that, if $0<\Delta\le \Delta_0(\epsilon)$
and if the condition
\begin{eqnarray}
\label{eq:condiz3}
\frac1{\beta^2}\ge a_2(\epsilon)\,\ell\left(\frac{a_2(\epsilon)|\log
F(0,\Delta)|}
{F(0,\Delta)
}\right)
\end{eqnarray}
is verified, then Eq. \eqref{eq:risulta1} holds.
\end{theorem}
\begin{rem}\rm
\label{r:56}
To be precise, in the statement of \cite[Th. 4]{cf:Ken} the condition
\eqref{eq:condiz3} is replaced by a different one (essentially, the
factor $|\log F(0,\Delta)|$ in the argument of $\ell(.)$ does not
appear). In this sense, the condition \eqref{eq:condiz3} under which
we prove here \eqref{eq:risulta1} is not the best possible one.  However,
for many ``reasonable'' and physically interesting choices of $L(.)$
in \eqref{eq:K}, Theorem \ref{th:a32} and Theorem 4 of \cite{cf:Ken}
are equivalent. In particular, if $\bP$ is the law of the returns to
zero of the simple random walk $\{S_n\}_{n\ge0}$ in one dimension,
i.e. $\tau=\{n\ge0:S_{2n}=0\}$, in which case $L(.)$ and $\tilde L(.)$
are asymptotically constant and $\ell(N)\sim a_3 \log N$, one sees
easily that \eqref{eq:condiz3} is verified as soon as
\begin{eqnarray}
\Delta\ge a_4(\epsilon)e^{-\frac{a_5(\epsilon)}{\beta^2}},
\end{eqnarray}
which is the same condition given in \cite{cf:Ken}.

Note, by the way, that in this case the difference $
h_c(\beta)-h^a_c(\beta)$ vanishes faster than any power of $
\beta$, for $\beta\searrow0$. This confirms the fact
that, even if the two critical points can be different, they cannot be
distinguished perturbatively.
\end{rem}

\subsection{Some open problems}
The results of previous section, while giving rigorous bases to
predictions based on the Harris criterion, leave various intriguing
gaps in our comprehension of the matter.  Let us list a few of them,
in random order:
\begin{itemize}
\item Let $\alpha<1/2$. Does there exist a $\beta_c<\infty$ such that
$h_c(\beta)\ne h_c^a(\beta)$ for $\beta>\beta_c$? If yes, how smooth is
$h_c(\beta)$ at $\beta_c$? Does $\nu$ equal $2-1/\alpha$ also for $\beta$
large?
\item Again, let $\alpha<1/2$ and look at Eq. \eqref{eq:svil}. Is it true that
$$
F(\beta,h^a_c(\beta)+\Delta)\ge F(0,\Delta)-\frac{\beta^2}2(\partial_\Delta
F(0,\Delta))^2(1+O(\beta^2))?
$$
\item Under the assumptions of Theorems \ref{th:a322} or \ref{th:a32},
does there exist positive values of $\beta$ for which quenched and
annealed critical points coincide? It is sort of reasonable to
conjecture that the answer is ``no'', at least for $\alpha>1/2$.
\end{itemize}
The reader might be tempted to think that such questions should be
easy to answer numerically. If so, he should have a look at
Ref. \cite{cf:CGG} where one gets an idea (in the context of random
heteropolymers at selective interfaces) of why numerical tests become
extremely hard in the neighborhood of the critical curve.

\begin{rem}\rm 
  Between the time these notes were written and the time they were
published, the above open problems have been to a
large extent solved.  In particular:
\begin{itemize}
\item in Ref. \cite{cf:T_aap} it was proven that 
for every $\alpha>0$, if $\gb$ is large enough and, say, $\go$ is Gaussian,
then $h_c(\gb)\ne h_c^a(\gb)$.
\item The question posed in open problem (2) has been answered
  positively in Ref. \cite{cf:GT_irrel}, although in a slightly weaker
sense.
\item In Ref. \cite{cf:DGLT} it was proven that as soon as $\alpha>1/2$ 
and $\gb>0$ one has $h_c(\gb)\ne h_c^a(\gb)$.
\end{itemize}
\end{rem}

\subsection{Proof of Theorems \ref{th:alpha<32}-\ref{th:a32}}
\label{sec:pfs}
We follow the approach of \cite{cf:T_qrc} which, with respect to that
of \cite{cf:Ken}, has the advantage of technical simplicity and of
being closely related to the interpolation ideas of Section
\ref{sec:UBq}. On the other hand, we  encourage the reader to look also
at the methods developed in \cite{cf:Ken}, which have the bonus of
extending in a natural way beyond the Gaussian case and of giving in
some cases sharper results (cf. Remarks \ref{r:54} and \ref{r:56}
above).

\medskip

A natural idea to show that quenched and annealed systems have
(approximately) the same free energy is to apply the {\sl second
moment method}: one computes $\bbE(Z_N(\beta,h))$ and $\bbE(
(Z_N(\beta,h))^2)$ and if it happens that the ratio
\begin{eqnarray}
\label{eq:secmom}
\frac{[\bbE Z_{N,\go}(\beta,h)]^2}{\bbE[ (Z_{N,\go}(\beta,h))^2]}
\end{eqnarray}
remains positive for $N\to\infty$, or at least it vanishes slower than
exponentially, it is not difficult to deduce that
$F(\beta,h)=F^a(\beta,h).$ This approach has turned out to be very
powerful for instance in controlling the high-temperature phase of the
Sherrington-Kirkpatrick mean field model in absence of magnetic field
\cite[Ch. 2.2]{cf:Tala_book}. However, this simple idea does not work
in our case and the ratio \eqref{eq:secmom} vanishes exponentially for
every $\beta,\Delta>0$. This is not surprising after all, since we
already know from Theorem \ref{th:rs} that quenched and annealed free
energy {\sl do not} coincide. There are two possible ways out of this
problem. One is to perform the second moment method not on the system
of size $N$ but on a smaller system whose size $N(\Delta)$ remains
finite as long as $\Delta$ is positive and fixed, and diverges only
for $\Delta\to0$. If $N(\Delta)$ is chosen to be the {\sl correlation
length of the annealed system}, one can see that on this scale the
ratio \eqref{eq:secmom} stays positive, so that
$F_{N(\Delta)}(\beta,h^a_c(\beta)+\Delta)\simeq
F_{N(\Delta)}(0,\Delta)$.  One is then left with the delicate problem
of glueing together many blocks of size $N(\Delta)$ to obtain an
estimate of the type $F(\beta,h^a_c(\beta)+\Delta)\ge
(1-\epsilon)F(0,\Delta)$ for the full free energy.  This is, in very
rough words, the approach of Ref. \cite{cf:Ken}.  The other possibility,
which we are going to present, is to abandon the second moment idea in
favor of a generalization of the {\sl replica coupling method}
\cite{cf:GuT_quadr} \cite{cf:T_qrc}.  This method was introduced in
\cite{cf:GuT_quadr} in the context of mean field spin glasses and
gives a very efficient control of the Sherrington-Kirkpatrick model at
high temperature ($\beta$ small), i.e., for weak disorder, which is
the same situation we are after here.

The two methods are in reality not orthogonal: they share the idea
that the important object to look at is the intersection of two
independent renewals $\tau^{(1)},\tau^{(2)}$. To see why this quantity
arises naturally, let us compute the second moment of the partition
function.  If $\tau^{(1)},\tau^{(2)}$ are independent renewal
processes with product law $\bP^{\otimes 2}(.)$, recalling the
definition $\Delta=h+\beta^2/2$, one can write
\begin{eqnarray}
  \bbE((Z_{N,\go}(\beta,h))^2)&=& \bbE \,\bE^{\otimes
2}\left(e^{\sum_{n=1}^N(\beta\go_n+h) ({\bf 1}_{\{n\in \tau^{(1)}\}}+
{\bf 1}_{\{n\in \tau^{(2)}\}})}{\bf 1}_{\{N\in \tau^{(1)}\}} {\bf
1}_{\{N\in \tau^{(2)}\}}\right)\\\nonumber &=&\bE^{\otimes
2}\left[e^{\Delta(|\tau^{(1)}\cap \{1,\ldots,N\}|+|\tau^{(2)}\cap
\{1,\ldots,N\}|)+\beta^2|\tau^{(1)} \cap\tau^{(2)}\cap
\{1,\ldots,N\}|}{\bf 1}_{\{N\in \tau^{(1)}\}} {\bf
1}_{\{N\in \tau^{(2)}\}}
\right].
\end{eqnarray}
Considering also that
$$
[\bbE\, Z_{N,\go}(\beta,h)]^2=
\bE^{\otimes 2}\left(e^{\Delta|(\tau^{(1)}\cap \{1,\ldots,N\}|+
|\tau^{(2)}\cap \{1,\ldots,N\}|)} {\bf 1}_{\{N\in \tau^{(1)}\}} {\bf
1}_{\{N\in \tau^{(2)}\}}\right) $$ one sees that the ratio
\eqref{eq:secmom} depends on the typical number of points that
$\tau^{(1)}$ and $\tau^{(2)}$ have in common up to time $N$.  One sees
also why this ratio has to vanish exponentially $N\to\infty$: as long
as $\Delta>0$ the renewals $\tau^{(i)}$, with law modified by the
factor $\exp(\Delta|\tau^{(1)}\cap \{1,\ldots,N\}|)$, are finitely
recurrent and therefore will have a number of intersections in
$\{1,\ldots,N\}$ which grows proportionally to $N$.

{\sl Proof of Theorem \ref{th:alpha<32}}.  The second inequality in
\eqref{eq:risulta1} is just Eq. \eqref{eq:jensen}.  As for the first
one, let $\Delta>0$ and recall identity \eqref{eq:id2}.  Define, in
analogy with \eqref{eq:Rt},
\begin{eqnarray}
  R_{N,\Delta}(t,\beta):=\frac1N\bbE\log
\med{e^{\sum_{n=1}^N(\beta\sqrt t\go_n
-t\beta^2/2)\delta_n}}_{\Delta,N}
\end{eqnarray}
for $0\le t\le 1$ (to the purpose of Theorem \ref{th:alpha<32} we do
not need the variational parameter $q$) where the measure
$\med._{N,\Delta}$ was defined after Eq. \eqref{eq:id1}.  Observe that
\begin{eqnarray}
\label{eq:obs1}
  R_{N,\Delta}(0,\beta)=0
\end{eqnarray}
while
\begin{eqnarray}
\label{eq:obs2}
  R_{N,\Delta}(1,\beta)=R_{N,\Delta}(\beta).
\end{eqnarray}
 As for the $t$-derivative one finds (just take \eqref{eq:derivq3} and
 put $q=0$):
\begin{eqnarray}
\label{eq:ipp1}
  \frac{\dd}{\dd t}R_{N,\Delta}(t,\beta)=-\frac{\beta^2}{2N}
\sum_{m=1}^N\bbE\left\{ \left(\frac{\med
{\delta_m\,e^{\sum_{n=1}^N(\beta\sqrt t\go_n-t\beta^2/2)\delta_n}}_{\Delta,N}}
{\med
{e^{\sum_{n=1}^N(\beta\sqrt t\go_n-t\beta^2/2)\delta_n}}_{\Delta,N}}\right)^2
\right\}.
\end{eqnarray}
Recall definition \eqref{eq:manag} (specialized to the case
$q=0$) of the random measure
$\med._{N,\Delta,t}$ and let
 $\med{.}_{N,\Delta,t}^{\otimes 2}$ be the product measure
acting on the pair $(\tau^{(1)},\tau^{(2)})$, while
$\delta^{(i)}_n:=\ind_{\{n\in\tau^{(i)}\}}$. Note that the two {\sl replicas}
$\tau^{(i)},i=1,2$ are subject to the {\sl same realization $\go$ of
disorder}. Then, one can rewrite
\begin{eqnarray}
\label{eq:derivata1}
 \frac{\dd}{\dd t}R_{N,\Delta}(t,\beta)=
-\frac{\beta^2}{2N}\bbE\sum_{m=1}^N
\med{\delta^{(1)}_m\delta^{(2)}_m}_{N,\Delta,t}^{\otimes 2}=-\frac{\beta^2}{2N}
\bbE\med{\left|\tau^{(1)}\cap \tau^{(2)}\cap \{1,\ldots,N\}\right|}
_{N,\Delta,t}^{\otimes 2} .
\end{eqnarray}
Since we need a lower bound for $R_{N,\Delta}(\beta)$ to prove
the first inequality in \eqref{eq:risulta1}, the fact that
this derivative is non-positive seems to go in the wrong direction.
Let us not lose faith and let us define, for $\lambda\ge0$,
\begin{eqnarray}
\label{eq:phi2}
 R^{(2)}_{N,\Delta}(t,\lambda,\beta)&:=&
\frac1{2N}\bbE \log
\med{e^{H_N(t,\lambda,\beta;\tau^{(1)},\tau^{(2)})}}_{N,\Delta}^{\otimes 2}
\\\nonumber
&:=&\frac1{2N}\bbE \log
\med{
e^{\sum_{n=1}^N(\beta\sqrt t \go_n-t\beta^2/2)(\delta^{(1)}_n+\delta^{(2)}_n)
+\lambda\beta^2\sum_{n=1}^N
\delta^{(1)}_n\delta^{(2)}_n}}_{N,\Delta}^{\otimes 2}
\end{eqnarray}
where the product measure $\med{.}_{N,\Delta}^{\otimes 2}$ acts on the
pair $(\tau^{(1)},\tau^{(2)})$.  The index ``$^{(2)}$'' refers to the
fact that this quantity involves two copies ({\sl replicas}) of the
system.  Observe that we are letting the two replicas interact through
a term which is positive, extensive (i.e., of order $N$) and closely
related to what appears in the right-hand side of Eq.
\eqref{eq:derivata1}.  Note also that
\begin{eqnarray}
   R^{(2)}_{N,\Delta}(0,\lambda,\beta)=
\frac1{2N}\log \med{
e^{\lambda\beta^2\sum_{n=1}^N
\delta^{(1)}_n\delta^{(2)}_n}}_{N,\Delta}^{\otimes 2},
\end{eqnarray}
while the factor $2$ in the denominator guarantees that
\begin{eqnarray}
\label{eq:310}
 R^{(2)}_{N,\Delta}(t,0,\beta)=R_{N,\Delta}(t,\beta).
\end{eqnarray}
Again via integration by parts (the computation is conceptually as
easy as the one which led to Eq. \eqref{eq:derivq3}),
\begin{eqnarray}
\label{eq:ddl}
\frac{\dd}{\dd t}R^{(2)}_{N,\Delta}(t,\lambda,\beta)&=&\frac{\beta^2}{2N}
\sum_{m=1}^N\bbE \frac{\med{\delta^{(1)}_m\delta^{(2)}_m
e^{H_N(t,\lambda,\beta;\tau^{(1)},\tau^{(2)})}}_{N,\Delta}^{\otimes
2} }{\med{e^{H_N(t,\lambda,\beta;\tau^{(1)},\tau^{(2)})}
}_{N,\Delta}^{\otimes 2}}\\\nonumber &&-
\frac{\beta^2}{4N}\sum_{m=1}^N\bbE
\left\{\left(\frac{\med{(\delta^{(1)}_m+\delta^{(2)}_m)
e^{H_N(t,\lambda,\beta;\tau^{(1)},\tau^{(2)})}}_{N,\Delta}^{\otimes
2} }{\med{e^{H_N(t,\lambda,\beta;\tau^{(1)},\tau^{(2)})}
}_{N,\Delta}^{\otimes 2}} \right)^2\right\}\\
\nonumber&&
\le \frac{\beta^2}{2N}\bbE
\sum_{m=1}^N\frac{\med{\delta^{(1)}_m\delta^{(2)}_m
e^{H_N(t,\lambda,\beta;\tau^{(1)},\tau^{(2)})}}_{N,\Delta}^{\otimes 2}
}{\med{e^{H_N(t,\lambda,\beta;\tau^{(1)},\tau^{(2)})}
}_{N,\Delta}^{\otimes 2}}= \frac{\dd}{\dd
\lambda}R^{(2)}_{N,\Delta}(t,\lambda,\beta).
\end{eqnarray}
This can be rewritten as
$$
\frac{\dd}{\dd t}R^{(2)}_{N,\Delta}(t,\lambda-t,\beta)\le 0
$$
which implies that, for every $0\le t\le 1$ and $\lambda$,
\begin{eqnarray}
\label{eq:daflusso}
R^{(2)}_{N,\Delta}(t,\lambda,\beta)\le R^{(2)}_{N,\Delta}(0,\lambda+t,\beta).
\end{eqnarray}
Going back to Eqs. \eqref{eq:ipp1} and the last equality in
\eqref{eq:ddl} and using the fact that for every convex function $\psi(.)$
one has $x\psi'(0)\le \psi(x)-\psi(0)$ one finds
\begin{eqnarray}
\label{eq:eta2}
&&   \frac{\dd}{\dd
t}\left(-R_{N,\Delta}(t,\beta)\right)=\frac{\dd}{\dd\lambda}
\left.R^{(2)}_{N,\Delta}(t,\lambda,\beta)\right|_{\lambda=0}
\le\frac{R^{(2)}_{N,\Delta}(t,2-t,\beta)-R^{(2)}_{N,\Delta}(t,0,\beta)}{2-t}.
\end{eqnarray}
Finally, using monotonicity of $R^{(2)}_{N,\Delta} (t,\lambda,\beta)$ with
respect to $\lambda$ and \eqref{eq:310}, one obtains the bound
\begin{eqnarray}
\frac{\dd}{\dd
t}\left(-R_{N,\Delta}(t,\beta)\right)
\le R^{(2)}_{N,\Delta}(0,2,\beta)+(-R_{N,\Delta}(t,\beta)),
\end{eqnarray}
where we used \eqref{eq:daflusso} and the fact that $2-t\ge1$ (of
course, we could have chosen $1+\eta-t$ instead of $2-t$ for some
$\eta>0$ in \eqref{eq:eta2} and the estimates would be modified in a
straightforward way).  We can now integrate with respect to $t$
between $0$ and $1$ this differential inequality (or use Gronwall's Lemma,
if you prefer) and, recalling
Eqs. \eqref{eq:obs2} and \eqref{eq:obs1}, we obtain
\begin{eqnarray}
\label{eq:integrating}
-(e-1)R^{(2)}_{N,\Delta}(0,2,\beta)\le R_{N,\Delta}(\beta)\le0.
\end{eqnarray}

Before we proceed, we would like to summarize what we did so far.  To
prove Theorem \ref{th:alpha<32} we need the {\sl lower bound}
$\lim_{N\to\infty} R_{N,\Delta}(\beta)\ge -\epsilon F(0,\Delta)$ but,
as in Section \ref{sec:UBq}, it seems that the interpolation method
gives rather {\sl upper bounds} on $R_{N,\Delta}(\beta)$. Then,
through the replica coupling trick we transferred this problem into
the problem of proving an {\sl upper bound} for a quantity,
$R^{(2)}_{N,\Delta}(t,\lambda,\beta)$, which is analogous to
$R_{N,\Delta}(\beta)$, except that it involves two interacting copies
of the system. Moreover, by throwing away a (complicated, but with a
definite sign) term in Eq. \eqref{eq:ddl}, we reduced to the problem
of bounding from above $R^{(2)}_{N,\Delta}(t,\lambda,\beta)$ computed
at $t=0$. In other words, we replaced the task of estimating from
below $R_{N,\Delta}(\beta)$ with that of estimating from above a
quantity which involves no quenched disorder, and which is therefore
easier to analyze. While this procedure might look a bit magic, the
basic underlying idea is the
following. $R^{(2)}_{N,\Delta}(t,\lambda,\beta)$ is obviously
non-decreasing as a function of $\lambda$. Suppose however that, for
some $\lambda>0$, $R^{(2)}_{N,\Delta}(t,\lambda,\beta)$ is not very
different from the value it has at $\lambda=0$ (of course, proving
this amounts to proving an {\sl upper} bound on
$R^{(2)}_{N,\Delta}(t,\lambda,\beta)$.)  Then, looking at the
definition \eqref{eq:phi2}, this means that the cardinality of the
intersection $\tau^{(1)}\cap \tau^{(2)}\cap
\{1,\ldots,N\}$ is typically not large and this, through
Eqs. \eqref{eq:obs1},
\eqref{eq:obs2} and \eqref{eq:derivata1} implies a {\sl lower} bound
on $R_{N,\Delta}(\beta)$.

Let us now restart from \eqref{eq:integrating} and note that
\begin{equation}
\begin{split}
\label{eq:split}
R^{(2)}_{N,\Delta}(0,2,\beta)\, &=\,
-F_N(0,\Delta)+
\frac1{2N}\log \bE^{\otimes 2}\left
(e^{2\beta^2\sum_{n=1}^N\delta^{(1)}_n\delta^{(2)}_n+\Delta
\sum_{n=1}^N(\delta^{(1)}_n+\delta^{(2)}_n)}\delta^{(1)}_N\delta^{(2)}_N
\right)
\\
&\le\,  -F_N(0,\Delta)+\frac{F_N(0,q\Delta)}q+
\frac1{2Np}\log \bE^{\otimes 2}
\left(e^{2p\beta^2\sum_{n=1}^N\delta^{(1)}_n\delta^{(2)}_n}
\right)
\end{split}
\end{equation}
where we used H\"older's inequality and the positive numbers $p$ and
$q$ (satisfying $1/p+1/q=1$) are to be determined. Taking the
thermodynamic limit,
\begin{multline}
\limsup_{N\to\infty}
  R^{(2)}_{N,\Delta}(0,2,\beta)\le \limsup_{N\to\infty}
\frac1{2Np}\log
\bE^{\otimes 2} \left(e^{2p\beta^2\sum_{n=1}^{N}\delta^{(1)}_n\delta^{(2)}_n}
\right)\\
+F(0,\Delta)\left(\frac1q \frac{F(0,q\Delta)}{F(0,\Delta)}-1\right).
\end{multline}
But we know from the expression  \eqref{eq:asint}
of the free energy of the homogeneous system and from the property
\eqref{eq:slow} of slow variation that, for every $q>0$,
\begin{eqnarray}
  \lim_{\Delta\to 0^+}\frac{F(0,q\Delta)}{F(0,\Delta)}=q^{1/\alpha}.
\end{eqnarray}
Therefore, taking $q=q(\epsilon)$ sufficiently close to (but strictly
larger than) $1$ and $\Delta_0(\epsilon)>0$ sufficiently small
one
has, uniformly on $\beta\ge0$ and on $0<\Delta\le \Delta_0(\epsilon)$,
\begin{eqnarray}
\label{eq:100}
 \limsup_{N\to\infty}R^{(2)}_{N,\Delta}(0,2,\beta)\le \frac{\epsilon}{e-1}
 F(0,\Delta)+
\limsup_{N\to\infty}\frac1{2Np(\epsilon)}\log
\bE^{\otimes 2}
\left(e^{2p(\epsilon)\beta^2\sum_{n=1}^{N}\delta^{(1)}_n\delta^{(2)}_n}
\right).
\end{eqnarray}
Of course, $p(\epsilon):=q(\epsilon)/(q(\epsilon)-1)<\infty$ as long
as $\epsilon>0$.  Note that, in view of \eqref{eq:integrating},
Theorem
\ref{th:alpha<32} would be proved if the second
term in the right-hand side of \eqref{eq:100} were zero. Up to now, we
have not used yet the assumption that $\alpha<1/2$ or that
\eqref{eq:orthat} holds, but now the right moment has come. The way this
assumption enters the game is that it guarantees that the renewal
$\tau^{(1)}\cap \tau^{(2)}$ is transient under the law $\bP^{\otimes
2}$.  Indeed,
\begin{eqnarray}
\label{eq:ricorrenza}
  \bE^{\otimes 2}
\left(\sum_{n\ge1}\ind_{n\in\tau^{(1)}\cap\tau^{(2)}}\right)=\sum_{n\ge1}
\bP(n\in\tau)^2<\infty
\end{eqnarray}
since, as proven in \cite{cf:doney},
\begin{eqnarray}
\label{eq:doney}
\bP(n\in\tau)\stackrel{n\to\infty}\sim\frac{C_\alpha}{L(n)n^{1-\alpha}}:=
\frac{\alpha\sin(\pi \alpha)}\pi
\frac1{L(n)n^{1-\alpha}}.
\end{eqnarray}
Actually, Eq. \eqref{eq:doney} holds more generally for $0<\alpha<1$ and
we will need it to prove Theorems \ref{th:a322} and \ref{th:a32}.

Transience and renewal properties of the process
of $\tau^{(1)}\cap \tau^{(2)}$ implies that
\begin{eqnarray}
   \label{eq:trans}
\bP^{\otimes 2}(|\tau^{(1)}\cap \tau^{(2)}|\ge k)\le (1-c)^k,
\end{eqnarray}
for some $0<c<1$: after each ``renewal epoch'', i.e., each 
point of $\tau^{(1)}\cap \tau^{(2)}$, the intersection renewal
has a positive probability $c$ of jumping to infinity.  Therefore,
there exists $\beta_1>0$ such that
\begin{eqnarray}
\label{eq:sup}
\sup_N \bE^{\otimes 2} \left(e^{2p(\epsilon)\beta^2\sum_{n=1}^{N}
\delta^{(1)}_n\delta^{(2)}_n}\right)<\infty
\end{eqnarray}
for every $\beta^2p(\epsilon)\le \beta_1^2$.  Together with
\eqref{eq:100} and \eqref{eq:id2}, this implies
\begin{eqnarray}
F(\beta,h^a_c(\beta)+\Delta)\ge (1-\epsilon)F(0,\Delta)
\end{eqnarray}
as soon as $\beta^2\le \beta_0^2(\epsilon):=\beta^2_1/p(\epsilon)$,
and therefore the validity of Theorem \ref{th:alpha<32}.
\hfill$\Box$

\medskip

{\sl Proof of Theorem \ref{th:a322}.} In what follows we assume that
$\Delta$ is sufficiently small so that $F(0,\Delta)\ll1$.  For
simplicity of exposition, we assume also that $L(.)$ tends to a
positive constant $L(\infty)$ at infinity (for the general case, which
is not significantly more difficult, cf. \cite{cf:T_qrc}).

If we try to repeat the proof of Theorem \ref{th:a322} in this case,
what goes wrong is that the intersection $\tau^{(1)}\cap\tau^{(2)}$ is
now recurrent, so that \eqref{eq:sup} does not hold any more. The
natural idea is then not to let $N$ tend to infinity at $\Delta$
fixed, but rather to work on a system of size $N(\Delta)$, which
diverges only when $\Delta\to0$, i.e., when the annealed critical
point is approached.  In particular, we let $N=N(\Delta):=c|\log
F(0,\Delta)|/F(0,\Delta)$ with $c>0$ large to be fixed later.  Note
also that this choice of $N(\Delta)$ is quite similar to that made in
\cite{cf:Ken}, where one applies the second moment method on a system
of size $c/F(0, \Delta)$ with $c$ large.  This choice has a clear
physical meaning: indeed, we will see in Section \ref{sec:xi} that the
correlation functions of the annealed system decay exponentially on
distances of order $1/F(0,\Delta)$ (the logarithmic factor in our
definition of $N(\Delta)$ should be seen as a technical necessity).

By the superadditivity property \eqref{eq:superadd} we have, in
analogy with \eqref{eq:id1},
 \begin{eqnarray}
 \label{eq:identity2}
   F(\beta,-\beta^2/2+\Delta)
 \ge F_{N(\Delta)}(0,\Delta)+R_{N(\Delta),\Delta}(\beta).
 \end{eqnarray}
To prove Theorem \ref{th:a322} we need to show that the first term in
the right-hand side of \eqref{eq:identity2} is essentially
$F(0,\Delta)$, while the second is not smaller than $-\epsilon
F(0,\Delta)$, in the range of parameters determined by condition
\eqref{eq:condiz}. The first fact is easy: as follows from Proposition
2.7 of \cite{cf:GT_ALEA}, there exists $a_6\in(0,\infty)$ (depending
only on the law $K(.)$ of the renewal) such that
  \begin{eqnarray}
\label{eq:alea}
  F_N(0,\Delta)\ge F(0,\Delta)-a_6\frac{\log N}N
  \end{eqnarray}
  for every $N$.    Choosing $c=c(\epsilon)$ large
 enough, Eq. \eqref{eq:alea} implies that
 \begin{eqnarray}
 \label{eq:fsc}
   F_{N(\Delta)}(0,\Delta)\ge (1-\epsilon)F(0,\Delta).
 \end{eqnarray}
 As for $R_{N(\Delta),\Delta}(\beta)$, we have from Eqs.
 \eqref{eq:integrating} and \eqref{eq:split}
 \begin{equation}
 \begin{split}
  \frac{R_{N(\Delta),\Delta}(\beta)}{e-1}\, \ge &\, -F(0,\Delta)
\left(\frac1q\frac{F(0,q\Delta)}{F(0,\Delta)}-1\right) -\epsilon
F(0,\Delta)\\ & \phantom{moveright} - \frac1{2N(\Delta)p}\log
\bE^{\otimes 2}
\left(e^{2p\beta^2\sum_{n=1}^{N(\Delta)} \delta^{(1)}_n\delta^{(2)}_n}
\right),
 \end{split}
 \end{equation}
where we used Eqs. \eqref{eq:fsc} and \eqref{eq:superadd} to bound
$-(1/q)F_{N(\Delta)}(0,q\Delta)+F_{N(\Delta)}(0,\Delta)$ from below.
Choosing again $q=q(\epsilon)$ we
obtain, for  $\Delta\le \Delta_0(\epsilon)$,
 \begin{eqnarray}
\label{eq:aus1}
\frac{R_{N(\Delta),\Delta}(\beta)}{e-1}\ge -2\epsilon F(0,\Delta) -
\frac1{2N(\Delta)p(\epsilon)}\log \bE^{\otimes 2}
\left(e^{2p(\epsilon)\beta^2\sum_{n=1}^{N(\Delta)}
\delta^{(1)}_n\delta^{(2)}_n} \right).
 \end{eqnarray} It was proven in \cite[Lemma 3]{cf:Ken} and
 \cite[Section 3.1]{cf:T_qrc} that if $1/2<\alpha<1$ there exists
 $a_{7}=\in(0,\infty)$, which depends in particular
on $L(\infty)$,  such that for every integers $N$
 and $k$ \begin{eqnarray} \label{eq:intersec1} \bP^{\otimes
 2}\left(\left|\tau^{(1)}\cap \tau^{(2)}\cap
\{1,\ldots,N\}\right|\ge k\right)\le
 \left(1-\frac{a_{7}}{N^{2\alpha-1}}\right)^k,
 \end{eqnarray}
which should be compared with \eqref{eq:trans}, valid for
$\alpha<1/2$.  Thanks to the geometric bound \eqref{eq:intersec1} we
have
\begin{eqnarray}
\label{eq:thanks}
  \bE^{\otimes 2} \left(e^{2p(\epsilon)\beta^2\sum_{n=1}^{N(\Delta)}
\delta^{(1)}_n\delta^{(2)}_n} \right)
&=&\sum_{k\ge0}\bP^{\otimes 2}\left(\sum_{n=1}^{N(\Delta)}
\delta^{(1)}_n\delta^{(2)}_n= k\right)e^{2p(\epsilon)\beta^2k}
\\\nonumber
&\le&
\left(
{1-e^{2\beta^2p(\epsilon)} \left(1-\frac{a_{7}}
{N(\Delta)^{2\alpha-1}}\right)}\right)^{-1},
\end{eqnarray}
whenever
$$
e^{2\beta^2p(\epsilon)}\left(1-\frac{a_{7}}{
N(\Delta)^{2\alpha-1}}\right)<1
$$
and this is of
course the case if
\begin{eqnarray}
\label{eq:stronger}
e^{2\beta^2p(\epsilon)}
\left(1-\frac{a_{7}}{ N(\Delta)^{2\alpha-1}}\right)\, \le\,
\left(1- \frac{a_{7}}{2N(\Delta)^{2\alpha-1}}\right).
\end{eqnarray}
At this point, using the definition of $N(\Delta)$ and
point (2) of Theorem \ref{th:hom}, it is not
difficult to see that there exists a positive constant
$a_{8}(\epsilon)$ such that
\eqref{eq:stronger} holds if
\begin{eqnarray}
\label{eq:horr}
\beta^2 p(\epsilon) \, &\le& \, a_{8}(\epsilon) \frac{
\Delta^{(2\alpha -1)/\alpha}}{
\left \vert \log F(0, \Delta) \right\vert ^{2\alpha -1}}.
\end{eqnarray}
Condition \eqref{eq:horr} is equivalent
to the first inequality in \eqref{eq:condiz}, for a suitable choice of
$a_1(\epsilon)$ and $\check L(.).$ As a consequence, for
$N(\Delta)$ sufficiently large (i.e., for $\Delta$ sufficiently small)
\begin{eqnarray}
\label{eq:asa}
   \frac1{2N(\Delta)p(\epsilon)}
\log \bE^{\otimes 2} \left(e^{2p(\epsilon)\beta^2\sum_{n=1}^{N(\Delta)}
 \delta^{(1)}_n\delta^{(2)}_n}
 \right)\le
\frac{F(0,\Delta)}{2c(\epsilon)p(\epsilon)|\log F(0,\Delta)|}
\log\left(\frac{2N(\Delta)^{2\alpha-1}}{a_{7}}
\right).
\end{eqnarray}

Recalling Eq. \eqref{eq:asint} one sees that, if $c(\epsilon)$ is
chosen large enough,
\begin{eqnarray}
\label{eq:together}
   \frac1{2N(\Delta)p(\epsilon)}
\log \bE^{\otimes 2} \left(e^{2p(\epsilon)\beta^2\sum_{n=1}^{N(\Delta)}
 \delta^{(1)}_n\delta^{(2)}_n}
 \right)\le \epsilon F(0,\Delta).
\end{eqnarray}
Together with Eqs. \eqref{eq:identity2},
\eqref{eq:fsc} and \eqref{eq:aus1}, this concludes the proof of the
theorem.
\hfill $\Box$

\medskip

{\sl Proof of Theorem \ref{th:a32}}. The proof is almost identical to
that of Theorem \ref{th:a322} and up to Eq. \eqref{eq:aus1} no changes
are needed. One has however to be careful with the geometric bound
\eqref{eq:intersec1}: in this case, it is not sufficient to replace
$\alpha$ by $1/2$, since the behavior at infinity of the slowly
varying function $L(.)$ in \eqref{eq:K} is here essential.  The
correct bound in this case is (cf. \cite[Lemma 3]{cf:Ken} and
\cite[Sec. 3.1]{cf:T_qrc}) \begin{eqnarray}
\label{eq:intersec2}
\bP^{\otimes
 2}\left(\sum_{n=1}^N \delta^{(1)}_n\delta^{(2)}_n\ge k\right)\le
 \left(1-\frac{a_{9}}{\ell(N)}\right)^k.  \end{eqnarray} for every
 $N$, for some $a_{9}>0$. We recall that $\ell(.)$ is the slowly
 varying function, diverging at infinity, defined by \eqref{eq:lN}.
 In analogy with Eq. \eqref{eq:thanks} one obtains then
\begin{eqnarray}
\label{eq:thanks2}
  \bE^{\otimes 2} \left(e^{2p(\epsilon)\beta^2\sum_{n=1}^{N(\Delta)}
\delta^{(1)}_n\delta^{(2)}_n} \right)\le
\left(1-e^{2\beta^2p(\epsilon)} \left(1-\frac{a_{9}}{\ell(N(\Delta))}\right)\right)^{-1}
\end{eqnarray}
whenever the right-hand side is positive.
Choosing $a_2(\epsilon)$ large enough one sees that if condition
\eqref{eq:condiz3} is fulfilled then
\begin{eqnarray}
e^{2\beta^2p(\epsilon)}\left(1-\frac{a_{9}}{\ell(N(\Delta))}\right)\le
\left(1-\frac{a_{9}}{2\ell(N(\Delta))}\right)
\end{eqnarray}
and, in analogy with \eqref{eq:asa},
\begin{eqnarray}
   \frac1{2N(\Delta)p(\epsilon)}
\log \bE^{\otimes 2} \left(e^{2(\epsilon)\beta^2\sum_{n=1}^{N(\Delta)}
 \delta^{(1)}_n\delta^{(2)}_n}
 \right)\le
\frac{F(0,\Delta)}{2c(\epsilon)p(\epsilon)|\log F(0,\Delta)|}
\log\left(\frac{2\ell(N(\Delta))}{a_{9}}\right).
\end{eqnarray}
From this estimate, for $c(\epsilon)$ sufficiently large one obtains
again \eqref{eq:together} and as a consequence the statement of
Theorem \ref{th:a32}.
\hfill $\Box$

\subsection{Smoothing effect of disorder (relevance for $\alpha>1/2$)}
\label{sec:a>}
Section \ref{sec:a<} was devoted to showing that, for
$\alpha<\alpha_c$, (weak) disorder is irrelevant, in that it does not
change the specific heat exponent $\nu$ and in that the transition
point coincides with the annealed one as long as $\beta$ is small.  We
saw also that for $\alpha_c\le\alpha<1$ quenched and annealed free
energies and critical points are very close (Theorems \ref{th:a322}
and \ref{th:a32}).  This might leave the reader with the doubt that
disorder might be irrelevant in this situation too. The purpose of the
present section is to show that this is not the case.

We start by recalling that via Theorem \ref{th:hom} and \eqref{eq:jensen}
we know that $F(\beta,h_c^a(\beta)+\Delta)\lesssim \Delta^{\max(1/\alpha,1)}$.
This bound is however quite poor: if we go back to \eqref{eq:RS} and
we choose $q=\Delta/\beta^2$ we obtain
\begin{eqnarray}
\label{eq:betterbut}
F(\beta,h^a_c(\beta)+\Delta)\le \frac{\Delta^2}{2\beta^2}
\end{eqnarray}
which is better, for $\Delta$ small and $\alpha>1/2$.
The point is however that, since one expects that $h^a_c(\beta)\ne
h_c(\beta)$ in this situation, \eqref{eq:betterbut} does not say anything
about the critical behavior of the quenched system; for this, we would need
rather an upper bound on $F(\beta,h_c(\beta)+\Delta)$. This is just the
content of the following result, which we state in the case of
Gaussian disorder:
\begin{theorem} \cite{cf:GT_smooth, cf:GT_prl}
For every $\beta>0$, $\alpha>0$ and $\Delta>0$ one has
\begin{eqnarray}
\label{eq:smooth}
F(\beta,h_c(\beta)+\Delta)\le \frac{(1+\alpha)}{2\beta^2}\Delta^2.
\end{eqnarray}
\label{th:smooth}
\end{theorem}
\begin{rem}\rm
Theorem \ref{th:smooth} actually holds beyond the Gaussian case; for
instance, in the case of bounded variables $\go_n$. In this case the
statement has to be modified in that the factor $2$ in that the
denominator in the right-hand side of \eqref{eq:smooth} is replaced by
$c:=c(\bbP)$, a constant which depends only on the disorder
distribution $\bbP$, and the results holds only provided $\Delta$ is
sufficiently small: $\Delta\le
\Delta_0(\bbP)$, see \cite{cf:GT_smooth}.
\end{rem}

\begin{rem}\rm
An obvious implication of Theorem \ref{th:smooth} is that $\nu\le 0$
as soon as $\beta>0$. In this sense, this result is much reminiscent
of what was proven in \cite{cf:CCFS,cf:CCFSlett} about the specific
heat exponent for the nearest-neighbor disordered Ising ferromagnet.

In particular, Theorem \ref{th:smooth} shows that the specific heat
exponent is modified by an arbitrary amount of disorder if
$\alpha>\alpha_c$: the phase transition is smoothed by randomness if
$\alpha>\alpha_c$ and becomes at least of second order (the effect is
particularly dramatic for $\alpha>1$, where the transition is of first
order for $\beta=0$).

It is also interesting to compare Theorem \ref{th:smooth} with the
celebrated result by M. Aizenman and J. Wehr \cite{cf:AW} which states
that first order phase transition in spin systems with discrete
spin-flip symmetry are smoothed by disorder as long as the spatial
dimension verifies $d\le 2$, while the same holds for $d\le 4$ if the
symmetry is continuous.

A less obvious consequence of Theorem \ref{th:smooth} is the following:
\begin{theorem}\cite{cf:T_jsp}\label{th:FSS}
Let $\beta>0$ and $0\le \alpha<\infty$. There exists $c>0$ such that
\begin{eqnarray}
\label{eq:2/3}
\lim_{N\to\infty}\bbE \,\bP_{N,\go}^{\beta,h_c(\beta)}\left
(|\tau\cap\{1,\ldots,N\}|\ge c N^{2/3}\log N\right)=0.
\end{eqnarray}
Moreover, under the assumptions of Theorem  \ref{th:alpha<32}, for
$\beta$ sufficiently small
\begin{eqnarray}
\label{eq:plus}
\lim_{N\to\infty}\bbE \,\bP_{N,\go}^{\beta,h_c(\beta)}\left
(|\tau\cap\{1,\ldots,N\}|\ge c N^{2\alpha/(1+\alpha)}\log N\right)=0.
\end{eqnarray}
\end{theorem}
This result should be read as follows. The fact that the transition is
at least of second order in presence of disorder implies already that
the Gibbs average of the contact fraction defined by \eqref{eq:cf}
tends to zero in the thermodynamic limit at the critical point. The
additional information provided by Theorem \ref{th:FSS} are finite-$N$
estimates on the size of $\tau\cap\{1,\ldots,N\}$ at criticality.
Whether the exponent $2/3$ in Eq. \eqref{eq:2/3} is optimal or not is an
intriguing open question.

Theorem \ref{th:FSS} was proven in \cite{cf:T_jsp}\footnote{ Theorem 3.1 in
\cite{cf:T_jsp} is formulated in the case of bounded random variables
$\go_n$, but it generalizes immediately to the Gaussian because the
basic ingredient one needs is the concentration inequality
\cite[Eq. (5.2)]{cf:T_jsp}, which holds in the case of Gaussian 
randomness as well.}
(together with more refined finite-size estimates on
$\bE_{N,\go}^{\beta,h}(|\tau\cap\{1,\ldots,N\}|)$ for $h-h_c(\beta)$
going to zero with $N$), apart from Eq. \eqref{eq:plus} which is a
consequence of
\cite[Th. 3.1]{cf:T_jsp} plus Theorem \ref{th:alpha<32} 
(cf. also Remark 3.2 in \cite{cf:T_jsp}).
\end{rem}

{\sl Proof of Theorem \ref{th:smooth} (sketch)} For a fully detailed
 proof we refer to \cite{cf:GT_smooth}.  In the case of
Gaussian disorder a simpler proof is hinted at in
\cite{cf:GT_prl} and fully developed in \cite[Section 5.4]{cf:GB}.

Here we give just a sketchy idea of why the transition cannot be of
first order when $\beta>0$. Assume by contradiction that
\begin{eqnarray}
 \label{eq:contrad}
F(\beta,h_c(\beta)+\Delta) \sim c \Delta\mbox{\;\; for\;\;} \Delta\to 0^+,
\end{eqnarray}
and consider the system at the critical point
$(\beta,h_c(\beta))$. Divide the system of size $N$ into $N/M$ blocks
$B_i$ of size $M$, with the idea that $1\ll M\ll N$. For a given
realization of $\go$ mark the blocks where the empirical average of
$\go$, i.e., $(1/M)\sum_{n\in B_i}\go_n$ equals approximately $
\Delta/\beta$. By standard large deviation estimates, there are typically
$\mathcal N_{marked}:=(N/M) e^{-M \Delta^2/(2\beta^2)}$ such blocks,
the typical distance between two successive ones being $D_{typ}:=M
e^{M \Delta^2/(2\beta^2)}$. It is a standard fact that if we take $M$
IID standard Gaussian variables and we condition on their empirical
average to be $\delta$, for $M$ large they are (roughly speaking)
distributed like IID Gaussian variables of variance $1$ {\sl and
average $\delta$}.  Therefore, in a marked block the system sees
effective thermodynamic parameters
$(\beta_{eff},h_{eff}):=(\beta,h_c(\beta)+\Delta)$.  Now we want to
show that the assumption \eqref{eq:contrad} leads to the (obviously
false) conclusion that $F(\beta,h_c(\beta))>0$. Indeed, let $\mathcal
S_\go$ be the set of $\tau$ configurations such that:
\begin{itemize}
\item there are no points of $\tau$ in unmarked blocks
\item the boundaries of all marked blocks belong to $\tau$.
\end{itemize}
Note that $\mathcal S_\go$ depends on disorder through the location
and the number of marked blocks, and that there is no restriction on
$\tau$ inside marked blocks.  One has the obvious bound
\begin{eqnarray}
\label{eq:gnaf}
F_N(\beta,h_c(\beta))\ge \frac1N \bbE\log
\bE\left(e^{\sum_{n=1}^N(\beta\go_n+h) \delta_n}{\bf 1}_{\{\tau\in
\mathcal S_\go\}}\delta_N\right).
\end{eqnarray}
But due to the definition of the set $\mathcal S_\go$, the restricted
free energy in the right-hand side of \eqref{eq:gnaf} gets (for $M$
large) a contribution $\mathcal N_{marked}\times (M/N)F(\beta,
h_c(\beta)+\Delta)$ from marked blocks, and an entropic term $\mathcal
N_{marked}/N\times
\log K(D_{typ})$ from the excursions between marked blocks.
Summing the two contributions, recalling the
asymptotic behavior \eqref{eq:K} of $K(.)$, the expression of
$\mathcal N_{marked}$ and $D_{typ}$ and taking the $N\to\infty$ limit
{\sl at $M$ large but fixed} one obtains then
\begin{eqnarray}
\label{eq:concl}
  F(\beta,h_c(\beta)+\Delta)\ge e^{-M \Delta^2/(2\beta^2)}\left(
F(\beta, h_c(\beta)+\Delta)-(1+\alpha)\frac{\Delta^2}{2\beta^2}
 \right).
\end{eqnarray}
Since the left-hand side of \eqref{eq:concl} is zero, for $\Delta$
small and $\beta>0$ this inequality is clearly in contradiction with
the assumption
\eqref{eq:contrad} that the transition if of first order (actually, even
with the assumption $F(\beta,h_c(\beta)+\Delta)\sim c \Delta^y$
with $y<2$).

\section{Correlation lengths and their critical behavior}
\label{sec:xi}
From certain points of view, the localized region $\cL$ is analogous
to the high-temperature phase of a spin system. Indeed, in this region
one can prove typical high-temperature results like the following:
free energy fluctuations are Gaussian on the scale $1/\sqrt N$
\cite{cf:AZ,cf:GT_ALEA}, the infinite-volume Gibbs measure is
almost-surely unique and ergodic \cite{cf:BidH}, the free energy is
infinitely differentiable, finite-size corrections to the infinite
volume free energy are of order $O(1/N)$, and truncated correlation
functions decay exponentially with distance \cite{cf:GT_ALEA}. In this
section we concentrate on the last point, which turns out to be more
subtle than expected, in particular when one approaches the critical
line.

In this section we assume that the random variables $\go_n$ are
bounded, because the results we mention have been proved in the
literature under this assumption. They should however reasonably
extend to more general situations, for instance to the Gaussian
case.

In the following, $\bP^{\beta,h}_{\infty,\go}(.)$ will denote the
infinite-volume Gibbs measure, defined as follows: first of all we
modify definitions \eqref{eq:Fo} and \eqref{eq:Boltz} replacing
$\sum_{n=1}^N(\beta\go_n+h)\delta_n$ by
$$
\sum_{n=-\lfloor N/2\rfloor
}^{\lfloor N/2\rfloor}(\beta\go_n+h)\delta_n, $$ where
$\{\go_n\}_{n\in \Z}$ are IID random variables, and then for a local
observable $f$, i.e., a function of $\tau$ which depends only on
$\tau\cap I$ with $I$ a finite subset of $\Z$, we let
\begin{eqnarray}
\label{eq:inftV}
  \bE^{\beta,h}_{\infty,\go}(f):=\lim_{N\to\infty}
\bE^{\beta,h}_{N,\go}(f).
\end{eqnarray}
Existence of the limit, in the localized phase, for almost every
disorder realization is proven in \cite{cf:GT_ALEA} (cf. also
\cite{cf:BidH}, where a DLR-like point of view is
adopted).\footnote{\label{nota:4} One might give a different
definition of the infinite-volume Gibbs measure, considering the
original system
\eqref{eq:Boltz} defined in $\{1,\ldots,N\}$ and taking a the
$N\to\infty$ limit of the average of local functions of $\tau\cap I$,
with $I$ a finite subset of $\N$.  In other words, with the first
procedure, Eq. \eqref{eq:inftV}, we are looking at the system in a
window which is situated in the bulk, very far away from both
boundaries. On the other hand, the second procedure is relevant if
one wants to study the system in the vicinity of one of the two boundaries
(and very far away from the other one).}

The definition of the correlation length $\xi$ contains always some
degree of arbitrariness, but conventional wisdom on universality
states that the critical properties of $\xi$, close to a second-order
phase transition, are insensitive to the precise definition.  There is
however a subtlety: in the case of disordered systems there
are two possible definitions of correlation lengths, which have no
reason to have the same critical behavior. Remaining for definiteness
in the framework of our disordered pinning models, one can first of
all define a {\sl (disorder-dependent) two-point function} as
\begin{eqnarray}
\label{eq:twop}
  \mathcal C_\go(k,\ell):=\bP_{\infty,\go}^{\beta,h}(k\in\tau|\ell
\in\tau)-
\bP^{\beta,h}_{\infty,\go}(k\in\tau).
\end{eqnarray}
In words, $\mathcal C_\go(k,\ell)$ quantifies how much the occurrence of
$\ell\in\tau$ influences the occurrence the event
$k\in\tau$.
It is then natural to define a correlation length $\xi$ as
\begin{eqnarray}
\label{eq:xity}
\frac1{\xi}:=  -\lim_{k\to\infty}\frac1k \log |\mathcal C_\go(k,0)|,
\end{eqnarray}
provided the limit exists. Note that $\xi$ depends on $(\beta,h)$ and,
in principle, on $\go$.
One can however define a different correlation length, $\xi^{av}$, as
\begin{eqnarray}
\label{eq:xiav}
  \frac1{\xi^{av}}:=-\lim_{k\to\infty}\frac1k \log \bbE|\mathcal C_\go(k,0)|.
\end{eqnarray}
In other words, $\xi$ (respectively, $\xi^{av}$) is the length over
which the two-point function (respectively, the averaged two-point
function) decays exponentially.  For simplicity, we will call $\xi$
the {\sl typical (or quenched) correlation length}, and $\xi^{av}$ the
{\sl average correlation length}, although it is important to keep in
mind that $\xi^{av}$ {\sl is not} the disorder-average of $\xi$
(indeed, in Section \ref{sec:xiI} we will see an example where $\xi$
is almost-surely constant but $\xi\ne \xi^{av}$).  It is interesting
that in the case of the one-dimensional quantum Ising chain with
random transverse field studied in \cite{cf:Fisher}, the two
correlation lengths are believed, on the basis of a renormalization
group analysis, to diverge at criticality {\sl with two different
critical exponents}.

A simple application of Jensen's inequality shows that $\xi^{av}\ge
\xi$.  This inequality can be interpreted on the basis of the
following intuitive argument. Divide all possible disorder
realizations into sets $A_m$ where the empirical average of $\go$ in
the region $\{1,\ldots,k\}$ is approximately $m$. Of course, for
$m\ne0$ $A_m$ is a large deviation-like event of probability $\simeq
\exp(-k m^2/2)$. Conditionally on $A_m$, the system sees a defect line
which is more attractive (if $m>0$) or more repulsive (if $m<0$) than
it should and therefore it is more localized (resp. more delocalized)
in this region than in the rest of the system. Therefore,
conditionally on $A_m$, we can expect that $\mathcal C_\go(k,0)$
behaves like $exp(-k/\xi(\beta,h+\beta m))$. In other words, we can
argue that (looking only at the exponential behavior)
\begin{eqnarray}
\label{eq:heu}
  \bbE \,\mathcal C_\go(k,0)\simeq \int \dd m\, e^{-k m^2/2}
e^{-k/ \xi(\beta,h+\beta m)}\simeq e^{k\max_{m}
\{-m^2/2-1/\xi(\beta,h+ \beta
m)\}}
\end{eqnarray}
for $k$ large.  Since $\xi$ should diverge when the critical point is
approached, it is reasonably decreasing in $h$ so that the value of
$m$ which realizes the maximum is strictly negative.  On the other
hand, when we take the limit without disorder average as in
\eqref{eq:xity}, the events $A_m$ with $m\ne0$ cannot contribute,
i.e.,  almost surely they do not occur for $k$ large enough, as follows from
the Borel-Cantelli lemma.

\subsection{Correlation length of the homogeneous model}
In the homogeneous case, $\beta=0$, the infinite-volume Gibbs measure
can be explicitly described (cf. \cite[Th. 2.3]{cf:GB}): under
$\bP^{0,h}_\infty(.)$, $\tau$ is a homogeneous\footnote{That is, its
law is invariant under translation on $\Z$. For instance,
$\bP^{0,h}_\infty(n,m\in\tau)= \bP^{0,h}_\infty(n+k,m+k\in\tau)$ for
every $k\in \Z$.}, positively recurrent (for $h>h_c(0)=0$) renewal on
$\Z$ such that
\begin{eqnarray}
\label{eq:tildek}
\bP^{0,h}_\infty(\inf\{k>0:k\in\tau\}=n|0\in\tau)=K(n)e^{-F(0,h)n}e^h
=:\tilde K_h(n)
\end{eqnarray}
and
$$ \bP^{0,h}_\infty(n\in\tau)=\frac1{\sum_{m\in  \N}m\tilde
K_h(m)}.
$$
Note that $\tilde K_h(.)$ is a probability on $\N$ (cf.
Eq. \eqref{eq:id3} and the discussion after it) with an exponential tail.
What we are interested in is the precise large-$n$ behavior of
$$\bP^{0,h}_\infty(n\in\tau|0\in\tau)
-
\frac1{\sum_{m\in  \N}m\tilde
K_h(m)},$$ 
i.e., a refinement of the renewal theorem (which simply
states that this quantity tends to zero for $n\to\infty$).

Let us for a moment widen our scope and consider a homogeneous,
positively recurrent renewal, with law  $\tilde\bP$, such
that the law of the distance between two successive points, 
denoted by $\tilde K(.)$, has exponential tail: say,
\begin{eqnarray}
\label{eq:ney}
  \lim_{n\to\infty} \frac1n\log \tilde K(n)=-z<0.
\end{eqnarray}
We {\sl do not} require for the moment that $\tilde K(.)$ is given by
\eqref{eq:tildek} with $K(.)$ in the class \eqref{eq:K}.  It is known
(cf. for instance \cite[Chapter VII.2]{cf:asmussen} and \cite{cf:Ney})
that, under condition \eqref{eq:ney}, there exist $r>0$ and $C<\infty$
such that
\begin{eqnarray}
  \label{eq:Kend2}
 \left|\tilde \bP(n\in \tau|0\in\tau)-\frac1{\sum_{m\in\N}m\tilde
K(m)}\right| \le C e^{-rn}.
\end{eqnarray}
However, the relation between $z$ and the largest possible $r$ in Eq.
\eqref{eq:Kend2}, call it $r_{max}$, is not known in general. A lot of
effort has been put by the queuing theory community in investigating
this point, and in various special cases it has been proven that
$r_{max}\ge z$ (see for instance \cite{cf:BL}, where power series
methods are employed and explicit upper bounds on the prefactor $C$
are given).  In even more special cases, for instance when $\tilde
\bP$ is the law of the return times to a particular state of a Markov
chain with some stochastic ordering properties, the optimal result
$r_{max}= z$ is proved (for details, see \cite{cf:LT,cf:T_jsp}, which
are based on coupling techniques).  However, the equality $r_{max}= z$
cannot be expected in general. In particular, if $\tilde K(.)$ is a
geometric distribution,
$$
\tilde K(n)=\frac{e^{-nc}}{e^c-1}
$$ with $c>0$, then one sees easily that the left-hand side of
\eqref{eq:Kend2} vanishes for every $n\in\N$ so that $r_{max}=\infty$,
while $z=c$. On the other hand, if for instance $\tilde K(1)=\tilde
K(2)=1/2$ and $\tilde K(n)=0$ for $n\ge3$, then $z=+\infty$ while $r$ is
finite.  These and other nice counter-examples are discussed in
\cite{cf:BL}.

In view of this situation, it is highly non-trivial that, restricting
to our original class of renewals, the following holds:
\begin{theorem}\cite{cf:DerrGiac}
\label{th:DerrGiac} Let $\tilde K_h(.)$ be given by \eqref{eq:tildek}
with $K(.)$ satisfying \eqref{eq:K} for some $\alpha>0$ and slowly
varying $L(.)$. Then, there exists $h_0>0$ such that, for every
$0<h<h_0$,
\begin{eqnarray}
  \label{eq:DerrGiac}
 \limsup_{n\to\infty}\frac1n\log 
\left|\bP^{0,h}_\infty(n\in \tau|0\in\tau)-
\frac1{\sum_{m\in\N}m\tilde K_h(m)}\right|
=-F(0,h)
\end{eqnarray}
and, more precisely,
\begin{eqnarray}
 \bP^{0,h}_\infty(n\in \tau|0\in\tau)-
\frac1{\sum_{m\in\N}m\tilde K_h(m)}\stackrel
{n\to\infty}\sim\frac{Q(n)e^{-F(0,h)n}}{4 [\sinh(h/2)]^2}
\end{eqnarray}
with $Q(.)$ such that $\sum_{j=1}^n Q(j)\stackrel
{n\to\infty}\sim L(n)/(\alpha n^\alpha)$.
\end{theorem}
It is important to emphasize that, even under assumption \eqref{eq:tildek},
this result would be false without the restriction of $h$ small.

In the light of \eqref{eq:DerrGiac}, it is quite natural to expect
(and in some case this can be proven, see Section \ref{sec:xiI}) that
in presence of disorder $\xi$ is still
proportional to the inverse of the free energy, at least
close to the critical point. But then, what about $\xi^{av}$?

\subsection{$\mu$ versus $F$}
To answer this question, we abandon for a while the correlation length and we
discuss the relation between free energy and another quantity which,
due to lack of a standard name, we will call simply $\mu$. This was
first introduced, to my knowledge, in \cite{cf:AZ} in the context of
random heteropolymers:
\begin{eqnarray}
  \label{eq:mu} \mu(\beta,h)=-\lim_{N\to\infty}\frac1N \log\bbE
  \left[\frac1 {Z_{N,\go}(\beta,h)}\right]
\end{eqnarray}
Existence of the limit in our context is easily proven by
super-additivity of $\log Z_{N,\go}(\beta,h)$ (see
\cite[Th. 2.5]{cf:GT_ALEA}).  An argument similar to \eqref{eq:F>0}
gives immediately $\mu\ge0$ while a simple application of Jensen's
inequality shows that $\mu(\beta,h)\le F(\beta,h)$. However, much more
than this is true:
\begin{theorem} \cite{cf:T_jsp}
\label{th:bounds}
For every $\beta>0$ there exists $0<c_3(\beta),c_4(\beta)<\infty$ such
that
  \begin{eqnarray}
\label{eq:bound_mu}
  0< c_3(\beta)\frac{F(\beta,h)^2}{\partial_h F(\beta,h)} <\mu(\beta,h)
<F(\beta,h)
  \end{eqnarray}
if $0<h-h_c(\beta)\le c_4(\beta)$.
\end{theorem}
In particular, the bounds in \eqref{eq:bound_mu} show that also $\mu$
vanishes continuously at the critical point, like the free energy.  If
we call $\eta_F$ and $\eta_\mu$ the critical exponents associated to
the vanishing of $F$ and $\mu$ for $h\to h_c(\beta)^+$, Theorem
\ref{th:bounds} implies the following bounds:
\begin{eqnarray}
(2\le)  \eta_F\le \eta_\mu\le \eta_F+1,
\end{eqnarray}
the inequality in parentheses being valid for $\beta>0$ thanks to
Theorem \ref{th:smooth}. Just to give a flavor of why $\mu$ is
relevant in the description of the system let us cite the following result.
Define first of all $\Delta_N$ as the largest gap
between points of $\tau$ in the system of length $N$:
  \begin{eqnarray}
    \Delta_N:=\max_{1\le i<j\le N}\{|i-j|:i\in\tau,j\in \tau,\{i+1,\ldots,
j-1\}\cap \tau=\emptyset\}.
  \end{eqnarray}
Then,
\begin{theorem}
\label{th:gap}
\cite{cf:GT_ALEA} Let $(\beta,h)\in\cL$. For every $\epsilon>0$,
  \begin{eqnarray}
\lim_{N\to\infty}    \bP_{N,\go}^{\beta,h}\left(\frac{1-\epsilon}{\mu(\beta,h)}
\le \frac{\Delta_N}{\log N}\le \frac{1+\epsilon}{\mu(\beta,h)}\right)=1
\mbox{\;\;in probability}.
  \end{eqnarray}
\end{theorem}

\subsection{Correlation lengths and free energy}
\label{sec:xiI}
To my knowledge, the only case where $\xi$ and $\xi^{av}$ can be
fully characterized even in presence of disorder is the one where
$K(.)$ is the law of the first return to zero of the one-dimensional
SRW conditioned to be non-negative. In other
words, let $\{S_n\}_{n=0,1,\ldots}$ be the SRW on $\Z$ started at
$S_0=0$ and let $\bP^{SRW}(.)$ denote its law.  We define
$K^{SRW,+}(n):=\bP^{SRW}(\inf\{k>0:S_k=0\}=2n|S_i\ge0\;\forall i)$. Go
back to Section \ref{rem:trans} for a motivation of this example as a
model of wetting of a $(1+1)$-dimensional substrate.  In this case, one
has the following
\begin{theorem}
\cite{cf:T_jsp}
\label{th:correlazioni}
Let $K(.)=K^{SRW,+}(.)$ and $\ell\in\Z$. For every $\beta\ge0$ and
$h>h_c(\beta)$,
\begin{eqnarray}
  \label{eq:corr_ave}
\frac1{\xi^{av}}=-\lim_{k\to\infty}\frac1k \log \bbE\, \mathcal
C_\go(\ell+k,k) =\mu(\beta,h)
\end{eqnarray}
and, $\bbP(\dd\go)$--a.s.,
\begin{eqnarray}
  \label{eq:corr_typ}
\frac1{\xi}=-\lim_{k\to\infty}\frac1k \log \mathcal C_\go(\ell+k,k)
=F(\beta,h).
\end{eqnarray}
\end{theorem}
With respect to Theorem \ref{th:DerrGiac}, this result is much less
sharp in that it catches only the exponential behavior of the
two-point function. However, note that in Theorem
\ref{th:correlazioni} $h-h_c(\beta)$ is not required to be small as in
Theorem \ref{th:DerrGiac}.  Note also that in
Eqs. \eqref{eq:corr_ave}, \eqref{eq:corr_typ} we have not taken the
absolute value of $C_\go(\ell+k,k)$: this is because, in this
particular case, one can prove that this quantity is non-negative
\cite{cf:T_jsp}.  Finally observe that, in view of
\eqref{eq:bound_mu}, the two correlation lengths are different. It
would be extremely interesting to know whether the two associated
critical exponents $\eta_F,\eta_\mu$
coincide or not.

\begin{rem}\rm
  Theorem \ref{th:correlazioni} does not coincide exactly with
\cite[Th. 3.5]{cf:T_jsp}, e.g., because in the latter
$\bP^{\beta,h}_{\infty,\go}(.)$ is the
infinite-volume Gibbs measure obtained from the system defined in
$\{1,\ldots,N\}$ letting $N\to\infty$ (cf. footnote
\ref{nota:4}). However, the proof of
\cite{cf:T_jsp} extends without difficulties to the result we stated
above. We remark also that the theorem holds as well in the case where
$K(n)=K^{SRW}(n):=\bP^{SRW}(\inf\{k>0:S_k=0\}=2n)$, i.e., the law of
the first return to zero of the unconditioned SRW. This follows from
the discussion in Section \ref{rem:trans} and from the fact that $K^{SRW}(n)=2
K^{SRW,+}(n)$.
\end{rem}

{\sl Proof of Theorem \ref{th:correlazioni} (sketch)}.  The proof of
Theorem
\ref{th:correlazioni} is based on a coupling argument. For simplicity
let $\bP^+(.):=\bP^{SRW}(.|S\ge0)$. One can then rewrite the two-point
function \eqref{eq:twop} as
\begin{eqnarray}
\nonumber
\mathcal C_\go(k,\ell)&=&\lim_{N\to\infty}\frac1{Z_{N,\go}(\beta,h)^{2}}
\\\nonumber
&&\times
\bE^{+,\otimes 2}\left[
e^{\sum_{n=-N/2}^{N/2}(\beta\go_n+h)({\bf 1}_{\{S^{(1)}_n=0\}}+
{\bf 1}_{\{S^{(2)}_n=0\}})}
\left({\bf 1}_{\{S^{(1)}_k=0\}}-{\bf 1}_{\{S^{(2)}_k=0\}}\right)
\vert S^{(1)}_\ell=0
\right],
\end{eqnarray}
where $S^{(1)},S^{(2)}$ are independent with law $\bP^+$.  Since the
 SRW conditioned to be non-negative is a Markov chain, the expectation
 in the right-hand side clearly vanishes if we condition on the event
 that there exists $\ell<i<k$ such that $S^{(1)}_i= S^{(2)}_i$. But
 (and here we use explicitly the condition $S_i\ge0$ and that two SRW
 trajectories which cross each other do necessarily intersect), if the
 complementary event happens then either $S^{(1)}$ or $S^{(2)}$ has no
 zeros in the interval $\{\ell+1,\ldots,k-1\}$.  As a consequence, one
 obtains
\begin{eqnarray}
  \bbE\, \mathcal C_\go(k,0)\le 2\bbE\,\bP^{\beta,h}_{\infty,\go}
(\tau\cap \{1,\ldots,k-1\}=\emptyset)
\end{eqnarray}
and it is not difficult to deduce from \eqref{eq:mu} that this
probability vanishes like $\exp(-k\mu(\beta,h))$ for
$k\to\infty$. For the opposite bound and for the proof of
\eqref{eq:corr_typ}  we refer to \cite{cf:T_jsp}. 
\medskip

In the general case where $\bP$ is not necessarily the law of the
returns of the SRW (or, in general, of any Markov chain), the
available results on correlation lengths in presence of disorder are
much less sharp and, above all, only correlation length {\sl upper
bounds} are known. At present, the best one can prove in general about
average correlation length is the following:
\begin{theorem}\cite{cf:T_xi}
\label{th:main_xi}
 Let $\epsilon>0$ and $(\beta,h)\in \mathcal L$. There exists
$C_1:=C_1(\epsilon,\beta,h)>0$ such that, for every $k\in\N$,
\begin{eqnarray}
  \label{eq:risultaAV}
\bbE \left|\mathcal C_\go(\ell+k,\ell)\right|\le
\frac1{C_1\mu(\beta,h)^{1/C_1}}
\exp\left(-k \,C_1\,\mu(\beta,h)^{1+\epsilon}\right).
\end{eqnarray}
The constant $C_1(\epsilon,\beta,h)$ does not vanish at the
critical line: for every bounded subset $B\subset \cL$ one has
$\inf_{(\beta,h)\in B}C_1(\epsilon,\beta,h)\ge C_1(B,\epsilon)>0$.
\end{theorem}
\begin{rem}\rm
\label{rimarca}
The necessity of introducing $\epsilon>0$ (i.e., of weakening the
upper bound with respect to the expected one)  is probably of
technical nature, as appears from the fact that for $\beta=0$ Theorem
\ref{th:main_xi} does not reproduce the sharp results
\eqref{eq:DerrGiac} which hold for the homogeneous case.

Observe that Theorem \ref{th:main_xi} is more than just an upper bound
on $\xi^{av}$. Indeed, thanks to the bound on the prefactor in front
of the exponential, Eq. \eqref{eq:risultaAV} says that the exponential
decay, with rate at least of order $\mu^{1+\epsilon}$, starts as soon
as $k\gg \mu^{-1-\epsilon}|\log \mu|$. This observation reinforces the
meaning of Eq. \eqref{eq:risultaAV} as an upper bound of order
$\mu^{-1}$ on the correlation length of disorder-averaged correlations
functions.
\end{rem}
About the typical correlation length the following can be proven:
\begin{theorem} \cite{cf:T_xi}
 \label{th:mainAS}
 Let $\epsilon>0$ and $(\beta,h)\in \mathcal L$. One has
 for every $k\in\N$
\begin{eqnarray}
\label{eq:risultaAS}
 \left|\mathcal C_\go(k,0)\right|\le
C_2(\go)
\exp\left(- k\,C_1\, F(\beta,h)^{1+\epsilon}\right),
\end{eqnarray}
where $C_1$ is as in Theorem \ref{th:main_xi}, while $C_2(\go):=
C_2(\go,\epsilon,\beta,h)$ is an almost surely finite random variable.
\end{theorem}

The proof of Theorems \ref{th:main_xi} and \ref{th:mainAS} relies on a
rather involved coupling/comparison argument. In simple (and
imprecise) words, one first approximates $K(.)$ with a new law $\tilde
K(.)$ which is the law of the returns to zero of a Markov process with
continuous trajectories (defined in terms of a Bessel process), and at
that point the coupling argument of last section can be applied. We
refer to \cite{cf:T_xi} for full details.

\appendix

\section{Two Tauberian results}

For completeness, we include without proof two Tauberian theorems
(i.e., results about the relation between the asymptotic behavior of a
function and of its Laplace transform) which we used in Section
\ref{sec:pfs}.  Given a function $Q:\N\to \R$, we define for $s\in\R$
$$
\hat Q(s):=\sum_{n\in\N}e^{-ns}Q(n)
$$
whenever the sum converges.

We begin with a (quite intuitive) fact:
\begin{theorem}\cite[Proposition 1.5.8]{cf:bingh}
\label{th:A1}
 If $\ell(.)$ is slowly varying and $\gamma>-1$ then
\begin{eqnarray}
\sum_{n=1}^N n^\gamma\ell(n)\stackrel{N\to\infty}\sim \frac{N^{\gamma+1}}
{\gamma+1}\ell(N).
\end{eqnarray}
\end{theorem}

Next we state Karamata's Tauberian theorem \cite[Th. 1.7.1]{cf:bingh}
which for our purposes may be formulated as follows:
\begin{theorem}
\label{th:A2}
Assume that $Q(n)\ge0 $ for every $n\in\N$, that $\ell(.)$ is slowly
varying and that $\rho\ge0$. The following are equivalent:
\begin{eqnarray}
\hat Q(s)\stackrel{s\searrow0}\sim \frac{\ell(1/s)}{s^\rho}
\end{eqnarray}
and
\begin{eqnarray}
\sum_{n=1}^N Q(n)\stackrel{N\to\infty}
\sim N^\rho \frac{\ell(N)}{\Gamma(1+\rho)}.
\end{eqnarray}
\end{theorem}
Recall that the function $\Gamma(z)$ can be defined, for $z>0$, as
$$
\Gamma(z)=\int_0^\infty  t^{z-1}e^{-t}\dd t.
$$

Finally, a theorem relating the Laplace transform of a law on the half-line
to its integrated tail (cf. \cite[Corollary 8.1.7]{cf:bingh}):
\begin{theorem}
\label{th:A3}
Let $X$ be an integer-valued random variables with law $\bP$ and
$Q(n):=\bP(X= n)$, $\ell(.)$ a slowly varying function and
$0\le\alpha<1$. The following are equivalent:
\begin{eqnarray}
1-\hat Q(s)\stackrel{s\searrow 0}\sim s^\alpha \ell(1/s)
\end{eqnarray}
and
\begin{eqnarray}
\bP(X>n)=\sum_{j>n}Q(j)\stackrel{n\to\infty}\sim
\frac{\ell(n)}{n^\alpha \Gamma(1-\alpha)}.
\end{eqnarray}
\end{theorem}

\section*{Acknowledgments}
I would like to thank Roman Koteck\'y for organizing the Prague Summer
School on Mathematical Statistical Mechanics and for inviting me to
give a course.  Learning and teaching there was an extremely
stimulating experience.

Some of the results described in these notes are based on joint work
with Giambattista Giacomin, to whom I am grateful for introducing me
to this subject, for countless motivating conversations, and also for
communicating to me the results of \cite{cf:DerrGiac} prior to
publication.

This work was supported in part by the GIP-ANR project JC05\_42461
({\sl POLINTBIO}) and my presence at the school was made possible
thanks to the support from the ESF-program ``Phase Transitions and
Fluctuation Phenomena for Random Dynamics in Spatially Extended
Systems''.

\end{document}